\documentclass{article}

\RequirePackage[OT1]{fontenc}
\RequirePackage{amsthm,amsmath,bm}
\RequirePackage[colorlinks,citecolor=blue,urlcolor=blue]{hyperref}

\usepackage{amsmath,amssymb, a4wide}
\usepackage{amssymb,latexsym}
\usepackage{dsfont, color}
\usepackage{verbatim}
\usepackage{cancel}
\usepackage{authblk}
\usepackage{graphicx}

\numberwithin{equation}{section}
\theoremstyle{plain}
\newtheorem{theorem}{Theorem}

\newtheorem{lemma}{Lemma}

\newtheorem{remark}{Remark}[section]
\newtheorem{example}{Example}

\newcommand{\R}{\mathbb{R}}

\def\argmin{\mathop{\mathrm{argmin}}}
\newcommand{\dd}{\mathrm{d}}
\newcommand{\vc}{\mathrm{vec}}
\newcommand{\vch}{\mathrm{vech}}

\newcommand{\bbeta}{\bm\beta}

\newcommand{\btheta}{\bm\theta}
\newcommand{\bxi}{\bm\xi}

\newcommand{\bmu}{\bm\mu}

\newcommand{\bSigma}{\mathbf{\Sigma}}

\newcommand{\by}{\mathbf{y}}
\newcommand{\bx}{\mathbf{x}}
\newcommand{\bs}{\mathbf{s}}
\newcommand{\bu}{\mathbf{u}}
\newcommand{\bv}{\mathbf{v}}

\newcommand{\bz}{\mathbf{z}}

\newcommand{\bb}{\mathbf{b}}

\newcommand{\bX}{\mathbf{X}}

\newcommand{\bV}{\mathbf{V}}

\newcommand{\bI}{\mathbf{I}}
\newcommand{\bA}{\mathbf{A}}
\newcommand{\bM}{\mathbf{M}}
\newcommand{\bC}{\mathbf{C}}
\newcommand{\bB}{\mathbf{B}}

\newcommand{\bL}{\mathbf{L}}
\newcommand{\bK}{\mathbf{K}}
\newcommand{\bO}{\mathbf{O}}

\newcommand{\bR}{\mathbf{R}}

\newcommand{\bN}{\mathbf{N}}
\newcommand{\bE}{\mathbf{E}}

\newcommand{\bJ}{\mathbf{J}}
\newcommand{\bT}{\mathbf{T}}

\newcommand{\E}{\mathbb{E}}

\title{Asymptotics of estimators for structured covariance matrices}
\author[1]{Hendrik Paul Lopuha\"a}
\affil[1]{\emph{Delft University of Technology}}
\date{\today}

\begin{document}
\maketitle
\begin{abstract}
We show that the limiting variance of a sequence of estimators for a structured covariance matrix
has a general form that appears as the variance of a scaled projection of a random matrix that is of radial type
and a similar result is obtained for the corresponding sequence of estimators for the vector of variance components.
These results are illustrated by the limiting behavior of estimators for a linear covariance structure
in a variety of multivariate statistical models.
We also derive a characterization for the influence function of corresponding functionals.
Furthermore, we derive the limiting distribution and influence function of scale invariant mappings
of such estimators and their corresponding functionals.
As a consequence, the asymptotic relative efficiency of different estimators for the shape component
of a structured covariance matrix can be compared by means of a single scalar
and the gross error sensitivity of the corresponding influence functions
can be compared by means of a single index.
Similar results are obtained for estimators of the normalized vector of variance components.
We apply our results to investigate how the efficiency, gross error sensitivity, and breakdown point
of S-estimators for the normalized variance components
are affected simultaneously by varying their cutoff value.
\end{abstract}

\section{Introduction}
\label{sec:introduction}
Covariance matrices describe the relationships and variability between different variables in a dataset.
When there is a known structure or pattern in these relationships,
structured covariance matrices can be estimated to capture and represent that structure.
The use of structured covariance matrices is a valuable tool
for modeling the underlying patterns and dependencies in multivariate data.
It provides a more nuanced understanding of the relationships between variables,
especially in scenarios where variables exhibit specific structures or patterns of correlation.
Structured covariance matrices are commonly used in the analysis of repeated measures,
longitudinal data, and multivariate data with a known underlying structure.
They are particularly useful when there are dependencies or correlations among different measurements or variables
and are widely used in various fields, including biology, medicine, psychology, and social sciences.

When a covariance matrix is unstructured and can be any positive definite symmetric matrix~$\bSigma$,
then the limiting behavior of covariance estimators $\bV_n$ for $\bSigma$ is well understood.
For example, if~$\bV_n$ is based on a sample $\by_1,\ldots,\by_n\in\R^k$ from a distribution with
an elliptically contoured density $|\bSigma|^{-1/2}g((\by-\bmu)^T\bSigma^{-1}(\by-\bmu))$,
then typically $\sqrt{n}(\bV_n-\bSigma)$ converges in distribution to a random matrix $\bN$
that has a  multivariate normal distribution with mean zero and variance
\begin{equation}
\label{eq:limiting variance unstructured}
\text{var}\{\text{vec}(\bN)\}
=
\sigma_1(\bI_{k^2}+\bK_{k,k})(\bSigma\otimes\bSigma)
+
\sigma_2\text{vec}(\bSigma)\text{vec}(\bSigma)^T,
\end{equation}
for some $\sigma_1\geq 0$ and $\sigma_2\geq -2\sigma_1/k$,
where $\otimes$ denotes the Kronecker product, $\bK_{k,k}$ is the commutation matrix,
and vec is the operator that stacks the columns of a matrix.
This form of limiting variance appears for many covariance estimators.
Tyler~\cite{tyler1982} gives several examples, including the sample covariance matrix,
and nicely explains that this general form will always appear
when $\bN$ is of radial type with respect to $\bSigma$.

The situation becomes different, when estimating a structured covariance matrix
$\bSigma=\bV(\btheta)$, where $\bV(\cdot)$ is a known covariance structure depending on
a vector $\btheta=(\theta_1,\ldots,\theta_\ell)$ of unknown variance components.
Asymptotic results for the maximum likelihood estimator of variance components
in linear models with Gaussian errors having a structured covariance matrix~$\bV(\btheta)$,
can be found in Hartley and Rao~\cite{hartley&rao1967},
Miller~\cite{miller1977}, and Mardia and Marshall~\cite{mardia-marshall1984}.
When scaled appropriately, the maximum likelihood estimator $\btheta_n$ is shown to be asymptotically normal
with mean~$\btheta$ and variance $\bJ^{-1}$, where
$\bJ_{ij}
=
\text{tr}(\bSigma^{-1}\bL_i\bSigma^{-1}\bL_j)/2$,
for $i,j,=1,\ldots,\ell$,
with $\bSigma=\bV(\btheta)$ and $\bL_i=\partial\bV(\btheta)/\partial\theta_i$.
By employing the vec-notation, the limiting covariance of~$\btheta_n$ can be expressed as
\[
2
\left(
\bL^T(\bSigma^{-1}\otimes\bSigma^{-1})\bL
\right)^{-1},
\]
where $\bL$ is the matrix with columns $\vc(\bL_1),\ldots,\vc(\bL_\ell)$.
According to the delta method the limiting covariance of $\vc(\bV(\btheta_n))$ is then given by
\[
2
\bL\left(
\bL^T(\bSigma^{-1}\otimes\bSigma^{-1})\bL
\right)^{-1}
\bL^T.
\]
Similar results have been obtained in Lopuha\"a \textit{et al}~\cite{lopuhaa-gares-ruizgazen2023} for the class of S-estimators
based on observations that follow a linear model with a structured covariance $\bSigma=\bV(\btheta)$, where~$\bV$ is
a linear function of~$\btheta$.
Under appropriate conditions, it holds that
$\sqrt{n}(\btheta_n-\btheta)$ is asymptotically normal with mean zero and variance
\begin{equation}
\label{eq:limiting variance theta structured}
2\sigma_1\Big(\bL^T\left(\bSigma^{-1}\otimes\bSigma^{-1}\right)\bL\Big)^{-1}
+
\sigma_2
\btheta\btheta^T,
\end{equation}
and $\sqrt{n}(\bV(\btheta_n)-\bSigma)$ converges in distribution to a random matrix $\bM$,
that has a  multivariate normal distribution with mean zero and variance
\begin{equation}
\label{eq:limiting variance structured}
\text{var}\{\text{vec}(\bM)\}
=
2\sigma_1
\bL\Big(\bL^T\left(\bSigma^{-1}\otimes\bSigma^{-1}\right)\bL\Big)^{-1}\bL^T
+
\sigma_2
\vc(\bSigma)\vc(\bSigma)^T.
\end{equation}
One of the objective of this paper is to show that this general form
will always appear when~$\bM$ is a scaled projection on the column space of $\bL$,
of a random matrix that is of radial type with respect to~$\bSigma$.
Moreover, we provide several examples of covariance estimators that exhibit this kind of limiting behavior.

Another objective concerns the asymptotic behavior of estimators for scale invariant mappings~$H$
of positive definite symmetric matrices.
For affine equivariant covariance estimators $\bV_n$ with asymptotic variance~\eqref{eq:limiting variance unstructured},
Tyler~\cite{tyler1983} shows that $H(\bV_n)$ has an asymptotic variance that only depends on the scalar $\sigma_1$.
When dealing with a structured covariance matrix, the covariance estimators are typically not affine equivariant
and have asymptotic variance~\eqref{eq:limiting variance structured}.
The second objective of this paper is to show that Tyler's result for affine equivariant covariance estimators,
remains true for estimators of a structured covariance matrix.
Moreover, we will establish a similar result for scale invariant mappings
$H(\btheta_n)$ of estimators for the vector of variance components.

An example of a scale invariant mapping is the shape component $\bV/|\bV|^{1/k}$.
A consequence of our results is that the asymptotic relative efficiency of estimators of the shape of
a structured covariance can be compared simply by comparing the
corresponding values for $\sigma_1$.
For affine equivariant covariance estimators,
this was already observed by Kent and Tyler~\cite{kent&tyler1996} and Salibi\'an \emph{et al}~\cite{SalibianBarrera-VanAelst-Willems2006}.
Similar properties will be shown to hold for the direction component $\btheta/\|\btheta\|$ corresponding to
the vector of variance components.

A final objective of this paper concerns the influence function of structured covariance functionals.
For affine equivariant covariance functionals, Croux and Haesbroeck~\cite{croux-haesbroeck2000} show that
the influence function at the multivariate normal is characterized by two real-valued functions.
Structured covariance functionals, however, are not necessarily affine equivariant.
We will show that such a characterization remains valid for structured covariance functionals
at any elliptically contoured distribution, 
and similarly for the variance components functional. 
A nice consequence is that the influence function of scale invariant mappings~$H$ of a
structured covariance functional~$\bV(\btheta(\cdot))$ or of $\btheta(\cdot)$ itself,
is characterized by a single real-valued function.
As such the gross-error-sensitivity (GES) is proportional to a single index,
which can be used to compare the GES of different shape functionals or different direction functionals.
Kent and Tyler~\cite{kent&tyler1996} already observed such a property for
the shape component of affine equivariant covariance functionals,
see also Salibi\'an \emph{et al}~\cite{SalibianBarrera-VanAelst-Willems2006}.

Except that our results have a merit of their own, they also enable the construction of MM-estimators
with auxiliary scale in linear mixed effects models and other linear models with structured covariances.
These estimators inherit the robustness of S-estimators considered in Lopuha\"a \textit{et al}~\cite{lopuhaa-gares-ruizgazen2023} and,
in contrast to the simpler version considered in Lopuha\"a~\cite{lopuhaa2023},
improve both the efficiency of the estimator of the fixed effects as well as the efficiency of the
estimator of the covariance shape component and of the direction of the vector of variance components.
Investigation of this version of MM-estimators will be postponed to a future manuscript, in which we will extend similar results that are already available
for unstructured covariances in the multivariate location-scale model, see Tatsuoka and Tyler~\cite{tatsuoka&tyler2000}
or Salibi\'an-Barrera \emph{et al}~\cite{SalibianBarrera-VanAelst-Willems2006},
and in the multivariate regression model, see Kudraszow and Maronna~\cite{kudraszow-maronna2011}.

The paper is organized as follows.
In Section~\ref{sec:projection} we show that the general forms of~\eqref{eq:limiting variance structured}
and~\eqref{eq:limiting variance theta structured} can be derived solely using a scaled projection of a random matrix that is of radial type.
In Section~\ref{sec:projection estimators} we investigate the limiting behavior of
estimators of a linear covariance structure in a variety of multivariate models.
We establish that these estimators asymptotically behave the same as a scaled projection
of a sequence of affine equivariant covariance estimators that are asymptotically of radial type.
In Section~\ref{sec:homogeneous mapping} we derive the limiting distribution of
scale invariant mappings of estimators of a linear covariance structure that are asymptotically normal,
and similarly for scale invariant mappings of estimators of the vector of variance components.
In Section~\ref{sec:IF} we derive a characterization for the influence function
of linearly structured covariance functionals and the corresponding functional
of variance components, and of scale invariant mappings thereof.
In Section~\ref{sec:application} we apply our results to investigate how
the efficiency, GES, and breakdown point of S-estimators
of the variance components are affected simultaneously,
when we vary the cut-off value of the rho-function that defines the S-estimator.
All proofs are postponed to an appendix at the end of the paper.

\section{Projection of a random matrix of radial type}
\label{sec:projection}
A random matrix $\bR$ is said to be of \textit{radial type},
if for any orthogonal matrix $\bO$,
the distribution of $\bO\bR\bO^T$ is the same as that of $\bR$.
The covariance structure of random matrices with a radial distribution
was first given by Mallows~\cite{mallows1961} in index form.
Tyler~\cite{tyler1982} gave the covariance structure in matrix form
and provided necessary conditions on its parameters.
A random matrix $\bN$ is said to be of radial type with respect to the positive definite symmetric matrix $\bSigma$,
if $\bSigma^{-1/2}\bN\bSigma^{-1/2}$ has a radial distribution.
If the first two moments of $\bN$ exist, then according to Corollary~1 in Tyler~\cite{tyler1982},
the variance of $\bN$  is given by~\eqref{eq:limiting variance unstructured}.

Consider a $k\times k$ structured covariance matrix $\bSigma=\bV(\btheta)$,
where $\bV$ is a known covariance structure that is a linear function of $\btheta=(\theta_1,\ldots,\theta_\ell)$,
a vector of unknown variance components.
Define
the $k^2\times \ell$ matrix
\begin{equation}
\label{def:L}
\bL
=
\left[
  \begin{array}{ccc}
 \vc(\bL_1)& \cdots & \vc(\bL_\ell) \\
  \end{array}
\right],
\quad
\bL_j=\partial\bV/\partial\theta_j,
\text{ for }
j=1,\ldots,\ell.
\end{equation}
Note that since $\bV$ is linear, we can write $\bSigma=\theta_1\bL_1+\cdots+\theta_\ell\bL_\ell$
and $\vc(\bSigma)=\bL\btheta$.
Furthermore, let $\Pi_L$ be the projection of a vector $\bx\in\R^{k^2}$ on the
column space of $\bL$, re-scaled by $\bSigma^{-1}\otimes\bSigma^{-1}$,
that is
\begin{equation}
\label{def:Projection}
\Pi_L\bx
=
\argmin_{\btheta\in\R^\ell}\,
(\bx-\bL\btheta)^T(\bSigma^{-1}\otimes\bSigma^{-1})(\bx-\bL\btheta).
\end{equation}
We then have the following theorem.
\begin{theorem}
\label{th:projection}
Let $\bN$ be a random matrix that is of radial type with respect to
a positive definite symmetric matrix $\bSigma$.
Suppose that $\bSigma=\bV(\btheta)$, for some $\btheta\in\R^{\ell}$, and that $\bV$ is linear
such that~$\bL$, as defined in~\eqref{def:L}, is of full column rank.
Let $\Pi_L$ be the projection defined in~\eqref{def:Projection}
and define the random matrix $\bM$ by $\vc(\bM)=\Pi_L\vc(\bN)$.
\begin{itemize}
\item[(i)]
If the first two moments of $\bN$ exist,
then there exist constants $\eta$, $\sigma_1$ and~$\sigma_2$ with $\sigma_1\geq 0$ and $\sigma_2\geq -2\sigma_1/k$, such that
$\E[\vc(\bM)]=\eta\vc(\bSigma)$ and
$\text{\rm var}(\vc(\bM))$ is given by~\eqref{eq:limiting variance structured}.
\item[(ii)]
If $\bT\in\R^\ell$ is the random vector, such that $\vc(\bM)=\bL\bT$, then
$\E[\bT]=\eta\btheta$ and
$\text{\rm var}(\bT)$ is given by~\eqref{eq:limiting variance theta structured}.
\end{itemize}
\end{theorem}
Note that the constants $\eta$, $\sigma_1$ and $\sigma_2$ have nothing to do with the projection $\Pi_L$,
but are inherited from the
variance~\eqref{eq:limiting variance unstructured} of the radial random matrix~$\bN$.
Their existence is guaranteed by Corollary~1 in Tyler~\cite{tyler1982}.

Examples of multivariate statistical models with a linear covariance structure
are linear mixed effects models.
But also linear models with errors generated by some autoregressive time series
may correspond to a linear covariance structure.
When $\bSigma$ is unstructured and can be any positive definite symmetric covariance matrix,
it can also be seen as a linear covariance structure $\bV(\btheta)$, where
$\btheta=\vch(\bSigma)$, with
\begin{equation}
\label{def:vech}
\vch(\bA)=(a_{11},\ldots,a_{k1},a_{22},\ldots,a_{kk}),
\end{equation}
is the unique $k(k+1)/2$-vector that stacks the columns of the lower triangle elements of a symmetric matrix $\bA$.
The matrix $\bL=\partial\vc(\bV)/\partial\btheta^T$ is then equal to the so-called duplication matrix~$\mathcal{D}_k$,
which is the unique $k^2\times k(k+1)/2$ matrix, with the properties
$\mathcal{D}_k\vch(\bA)=\vc(\bA)$ and $(\mathcal{D}_k^T\mathcal{D}_k)^{-1}\mathcal{D}_k^T\vc(\bA)=\vch(\bA)$.
Moreover, from the properties of $\mathcal{D}_k$ (e.g., see
Magnus and Neudecker~\cite[Ch.~3, Sec.~8]{magnus&neudecker1988}), it follows that
\begin{equation}
\label{eq:property duplication matrix}
\mathcal{D}_k\left(
\mathcal{D}_k^T
\left(\bSigma^{-1}\otimes\bSigma^{-1}\right)
\mathcal{D}_k
\right)^{-1}
\mathcal{D}_k^T
=
\frac12\left(
\bI_{k^2}+\bK_{k,k}
\right)
\left(\bSigma\otimes\bSigma\right).
\end{equation}
In this case, the expression~\eqref{eq:limiting variance structured}
with $\bL=\mathcal{D}_k$
coincides with the expression~\eqref{eq:limiting variance unstructured}.

\section{Projections of estimators of radial type}
\label{sec:projection estimators}
A sequence $\{\bN_n\}$ of $k\times k$ symmetric estimators for $\bSigma$ is said to be
\emph{asymptotically of radial type}
if there exists a sequence of real numbers $a_n$ increasing to infinity,
such that $a_n(\bN_n-\bSigma)\to\bN$ in distribution
with $\bN$ being of radial type with respect to $\bSigma$,
see Tyler~\cite{tyler1982}.
In a large class of multivariate statistical models,
for estimators $\bV_n$ of a linearly structured covariance matrix,
it turns out that the limiting behavior of $\vc(\bV_n)$ is the same as that of the projection $\Pi_L\vc(\bN_n)$
of a random matrix $\bN_n$ that is asymptotically of radial type with respect to $\bSigma$,
where $\Pi_L$ is defined in~\eqref{def:Projection}.
We illustrate this behavior in the following linear model with a structured covariance.

Consider independent observations $\bs_1,\ldots,\bs_n\in\R^k\times\R^{kq}$
with distribution $P$,
where $\bs_i=(\by_i,\bX_i)$, $i=1,\ldots,n$,
for which we assume the following model
\begin{equation}
\label{def:linear model}
\by_i
=
\bX_i\bbeta+\bu_i,
\quad
i=1,\ldots,n,
\end{equation}
where $\by_i\in\R^{k}$,
$\bbeta\in\R^q$ is an unknown parameter vector,
$\bX_i\in\R^{k\times q}$ is a known design matrix, and
$\mathbf{u}_i\in\R^{k}$ are unobservable independent mean zero random vectors with
covariance matrix $\bV\in\text{PDS}(k)$,
the class of positive definite symmetric $k\times k$ matrices.
Suppose that the distribution~$P$ for random variable $\bs=(\by,\bX)$ is such that~$\by\mid\bX$ has an elliptically contoured density
\begin{equation}
\label{eq:elliptical}
f_{\bmu,\bSigma}(\by)
=
|\bSigma|^{-1/2}
g\left(
(\by-\bmu)^T\bSigma^{-1}(\by-\bmu)
\right),
\end{equation}
where $\bmu=\bX\bbeta_0$ and $\bSigma=\bV(\btheta_0)=\theta_{01}\bL_1+\cdots+\theta_{0\ell}\bL_\ell$,
for some vector~$\btheta_0\in\R^\ell$ of variance components.
This setup includes several multivariate statistical models of interest.
One possibility is the linear mixed effects model, in which the random effects together with the measurement error
yields a specific covariance structure.
Other covariance structures may arise, for example if the~$\bu_i$ are the outcome of a time series.
Note that this setup also allows models with an unstructured covariance matrix,
such as the multivariate location-scale model or the multivariate regression model.
See e.g., Jennrich and Schluchter~\cite{jennrich&schluchter1986} or
Fitzmaurice \emph{et al}~\cite{fitzmaurice-laird-ware2011}, for different possible covariance structures,
and Lopuha\"a \emph{et al}~\cite{lopuhaa-gares-ruizgazen2023}, who provide a uniform treatment of S-estimators in these models.

Estimators $\bxi_n=(\bbeta_n,\btheta_n)$ for $\bxi_0=(\bbeta_0,\btheta_0)$ are typically solutions of estimating equations of the following type
\begin{equation}
\label{eq:score equations M}
\int \Psi(\bs,\bxi)\,\text{d}\mathbb{P}_n(\bs)=\mathbf{0},
\end{equation}
where $\mathbb{P}_n$ denotes the empirical measure corresponding to
$\bs_1,\ldots,\bs_n$, and where $\Psi=(\Psi_{\bbeta},\Psi_{\btheta})$, with
\begin{equation}
\label{eq:Psi function}
\begin{split}
\Psi_{\bbeta}(\bs,\bxi)
&=
w_1(d)\bX^T\bV^{-1}(\by-\bX\bbeta)\\
\Psi_{\btheta}(\bs,\bxi)
&=
\bL^T(\bV^{-1}\otimes\bV^{-1})
\vc
\left\{
w_2(d)
(\by-\bX\bbeta)(\by-\bX\bbeta)^T
-
w_3(d)
\bV
\right\},
\end{split}
\end{equation}
where $d^2=(\by-\bX\bbeta)^T\bV^{-1}(\by-\bX\bbeta)$,
and where we write $\bV$ for $\bV(\btheta)$.
We give some examples below.
Furthermore, typically $\bxi_n$ will then converge to a solution of
the corresponding population equation
\begin{equation}
\label{eq:score equations M for P}
\int \Psi(\bs,\bxi)\,\text{d}P(\bs)=\mathbf{0}.
\end{equation}
Let $\bV_n=\bV(\btheta_n)$.
From the estimating equations~\eqref{eq:score equations M} for $\bxi_n$,
we will establish that~$\vc(\bV_n)$ is asymptotically equivalent with $\Pi_L\vc(\bN_n)$,
for some $\bN_n$ that is asymptotically of radial type
and~$\Pi_L$ defined in~\eqref{def:Projection}.
To this end, we require the following conditions
\begin{quote}
\begin{itemize}
\item[(C1)]
$w_i(s)$ is of bounded variation and continuously differentiable, for $i=1,2,3$;
\item[(C2)]
$w_1'(s)s^2$, $w_2'(s)s^3$, and $w_3'(s)s^2$ are bounded;
\item[(C3)]
$\E_{\mathbf{0},\bI_k}\Big[w_2'(\|\bz\|)\|\bz\|^3+k(k+2)w_3(\|\bz\|)\Big]\ne0$ and
$\E_{\mathbf{0},\bI_k}\Big[w_2'(\|\bz\|)\|\bz\|^3+2kw_3(\|\bz\|)-kw_3'(\|\bz\|)\|\bz\|\Big]\ne0$,
\end{itemize}
\end{quote}
where $\E_{\mathbf{0},\bI_k}$ denotes the expectation with respect to
density~\eqref{eq:elliptical} with parameters $(\bmu,\bSigma)=(\mathbf{0},\bI_k)$.
Condition (C3) is to ensure the existence of the scalars $\sigma_1$ and $\sigma_2$ in Theorem~\ref{th:expansion}.
Maronna~\cite{maronna1976} and Tyler~\cite{tyler1982} consider M-estimators for multivariate
location and covariance.
Estimating equations for these estimators would correspond to $\Psi_{\btheta}$ without
the factor $\bL^T(\bV^{-1}\otimes\bV^{-1})$ (see Example~\ref{example:MLE} below) and $w_3=1$.
Moreover, they assume that $w_2'$ is non-negative, which obviously implies~(C3).

\begin{theorem}
\label{th:expansion}
Let $P$ be a distribution for random variable $\bs=(\by,\bX)$, such that
$\by\mid\bX$ has an elliptically contoured density~\eqref{eq:elliptical},
with parameters $\bmu=\bX\bbeta_0$ and $\bSigma=\bV(\btheta_0)$, for a linear covariance structure~$\bV$.
Let~$\bxi_n$ and $\bxi_0$ be solutions of~\eqref{eq:score equations M}
and~\eqref{eq:score equations M for P}, respectively, and suppose that $\bxi_n\to\bxi_0$ in probability.
Suppose that~$\E\|\bs\|^4<\infty$ and that~$\bX$ has full rank
with probability one.
If $w_1$, $w_2$, and $w_3$ satisfy (C1)-(C3), then
there exists a sequence $\{\bN_n\}$ of random matrices, such that
\[
\sqrt{n}
\left\{
\vc(\bV_n)-\vc(\bSigma)
\right\}
=
-
\Pi_L\vc\left\{
\sqrt{n}(\bN_n-\E[\bN_n])
\right\}
+
o_P(1),
\]
where $\Pi_L$ is defined in~\eqref{def:Projection}.
Moreover, $\sqrt{n}(\bN_n-\E[\bN_n])\to\bN$ in distribution, where $\bN$ is
a random matrix that has a multivariate normal distribution with mean zero and
variance~\eqref{eq:limiting variance unstructured}, with
\[
\begin{split}
\sigma_1
&=
\frac{k(k+2)\E_{\mathbf{0},\bI_k}\left[w_2(\|\bz\|)^2\|\bz\|^4\right]}{
\Big(\E_{\mathbf{0},\bI_k}
\Big[
w_2'(\|\bz\|)\|\bz\|^3+k(k+2)w_3(\|\bz\|)
\Big]\Big)^2}\\
\sigma_2
&=
-\frac{2}{k}\sigma_1
+
\frac{4\E_{\mathbf{0},\bI_k}\left[\Big(w_2(\|\bz\|)\|\bz\|^2-kw_3(\|\bz\|)\Big)^2\right]}{
\left(\E_{\mathbf{0},\bI_k}\Big[w_2'(\|\bz\|)\|\bz\|^3+2kw_3(\|\bz\|)-kw_3'(\|\bz\|)\|\bz\|\Big]\right)^2
}.
\end{split}
\]
\end{theorem}

\begin{remark}
\label{rem:explicit expression N}
From the proof of Theorem~\ref{th:expansion} one can obtain the following explicit expression for~$\bN_n$:
\[
\bN_n=\frac{1}{n}\sum_{i=1}^n
\left\{
v_1(d_i)(\by_i-\bX_i\bbeta_0)(\by_i-\bX_i\bbeta_0)^T
-
v_2(d_i)\bSigma
\right\},
\]
where $d_i^2=(\by_i-\bX_i\bbeta_0)^T\bSigma^{-1}(\by_i-\bX_i\bbeta_0)$,
and
\[
\begin{split}
v_1(s)
&=
\frac{w_2(s)}{\gamma_1};\\
v_2(s)
&=
\frac{-\gamma_2w_2(s)s^2+\gamma_1w_3(s)}{\gamma_1(\gamma_1-k\gamma_2)},
\end{split}
\]
where $\gamma_1$ and $\gamma_2$ are defined in~\eqref{def:gamma12}.
Note that $\gamma_1$ and $\gamma_1-k\gamma_2$ are precisely the quantities that appear in condition~(C3).
\end{remark}

The random matrix $\bN$ in Theorem~\ref{th:expansion} is of radial type with respect to $\bSigma$.
This follows from the fact that $\bR=\bSigma^{-1/2}\bN\bSigma^{-1/2}$ is multivariate normal with mean zero and variance
\[
\begin{split}
\text{var}\{\text{vec}(\bR)\}
&=
(\bSigma^{-1/2}\otimes\bSigma^{-1/2})
\text{var}
\left\{
\text{vec}(\bN)
\right\}
(\bSigma^{-1/2}\otimes\bSigma^{-1/2})\\
&=
\sigma_1(\bI_{k^2}+\bK_{k,k})
+
\sigma_2\text{vec}(\bI_k)\text{vec}(\bI_k)^T.
\end{split}
\]
This immediately gives that for any orthogonal matrix $\bO$,
the matrix $\bO\bR\bO^T$ is multivariate normal with mean zero and the same variance.
From Theorem~\ref{th:expansion} it follows that $\sqrt{n}\left\{\vc(\bV_n)-\vc(\bSigma)\right\}$ is
asymptotically normal with mean zero and a variance that is the same as the variance of $\vc(\bM)=\Pi_L\vc(\bN)$.
According to Theorem~\ref{th:projection} this variance is of the type given by~\eqref{eq:limiting variance structured}.
Furthermore, if we write $\vc(\bM)=\bL\bT$, then
\[
\sqrt{n}(\btheta_n-\btheta_0)
=
(\bL^T\bL)^{-1}\bL^T\sqrt{n}\left\{\vc(\bV_n)-\vc(\bSigma)\right\}
\to
\bT,
\]
in distribution,
where $\bT$ is multivariate normal with mean zero and variance given by~\eqref{eq:limiting variance theta structured}.

\subsection{Examples}
\label{subsec:examples}
We discuss some examples of multivariate statistical models that are covered by the setup in~\eqref{def:linear model},
in which the estimators $(\bbeta_n,\btheta_n)$ are solutions of estimating equation~\eqref{eq:score equations M}
for particular functions~$w_1$, $w_2$, and $w_3$.
In the Appendix we provide a detailed derivation $\sigma_1$ and $\sigma_2$ for specific special cases
and show that their expressions coincide with the ones in Tyler~\cite{tyler1982} and Lopuha\"a \emph{et al}~\cite{lopuhaa2023}.

\begin{example}[Maximum likelihood for multivariate normal]
\label{example:LS}
Suppose that $(\by_1,\bX_1),\ldots,(\by_n,\bX_n)$ are independent,
such that $\by_i\mid\bX_i\sim N_k(\bX_i\bbeta_0,\bV(\btheta_0))$.
The loglikelihood is then given by
\[
\mathcal{L}=
-\frac{nk}{2}\log(2\pi)
-\frac{n}2\log|\bV(\btheta)|
-\frac12
\sum_{i=1}^{n}
(\by_i-\bX_i\bbeta)^T\bV(\btheta)^{-1}(\by_i-\bX_i\bbeta).
\]
Setting the partial derivatives $\partial \mathcal{L}/\partial\bbeta$
and $\partial \mathcal{L}/\partial\theta_j$ equal to zero gives the following estimating equations
\begin{equation}
\label{eq:score equations LS}
\begin{split}
\frac1n
\sum_{i=1}^{n}
\bX_i^T\bV^{-1}(\by_i-\bX_i\bbeta)
&=\mathbf{0},\\
\frac1n
\sum_{i=1}^{n}
\left\{
(\by_i-\bX_i\bbeta)^T
\bV^{-1}\bL_j\bV^{-1}(\by_i-\bX_i\bbeta)
-
\text{\rm tr}(\bV^{-1}\bL_j)
\right\}&=0,
\end{split}
\end{equation}
for $j=1,\ldots,\ell$,
where we write~$\bV$ for~$\bV(\btheta)$.
By using the vec-notation and $\bL$ as defined in~\eqref{def:L},
we can combine the partial derivatives with respect to~$\theta_j$ in the second line of~\eqref{eq:score equations LS}
as follows
\begin{equation}
\label{eq:score equations LS vec}
\bL^T(\bV^{-1}\otimes\bV^{-1})
\vc
\left\{
\frac1n
\sum_{i=1}^{n}
(\by_i-\bX_i\bbeta)(\by_i-\bX_i\bbeta)^T
-
\bV
\right\}
=
\mathbf{0}.
\end{equation}
It follows that the maximum likelihood estimator $(\bbeta_n,\btheta_n)$
satisfies~\eqref{eq:score equations M} and $(\bbeta_0,\btheta_0)$ satisfies~\eqref{eq:score equations M for P}, where~$\Psi$
is defined in~\eqref{eq:Psi function} with $w_1(s)=w_2(s)=w_3(s)=1$.
Theorem~\ref{th:expansion} applies and one finds~$\sigma_1=1$ and $\sigma_2=0$.

When each $\bX_i=\bI_k$, for $i=1,\ldots,n$, then the model~\eqref{def:linear model} reduces to the multivariate
location-scale model.
If $\bSigma$ is unstructured, then $\bSigma=\bV(\btheta_0)$,
with $\btheta_0=\vch(\bSigma)$ and $\bL=\partial\vc(\bV(\btheta_0))/\partial\btheta^T$ is equal to the duplication matrix~$\mathcal{D}_k$.
In this case, we can remove the factor $\bL^T(\bV^{-1}\otimes\bV^{-1})$ from~\eqref{eq:score equations LS vec},
and $\bV_n$ is simply the sample covariance of~$\by_1,\ldots,\by_n$.
This example then coincides with Example~1 in Tyler~\cite{tyler1982}.
\end{example}

\begin{example}[M-estimators]
\label{example:MLE}
As mentioned in Example~\ref{example:LS},
when each $\bX_i=\bI_k$, for $i=1,\ldots,n$, and~$\bSigma$ is unstructured,
then the model~\eqref{def:linear model} reduces to the multivariate location-scale model
and we can remove the factor $\bL^T(\bV^{-1}\otimes\bV^{-1})$ from $\Psi_{\btheta}$ in~\eqref{eq:score equations M}.
In that case, estimating equations~\eqref{eq:score equations M} are equivalent to equations (1.1)-(1.2) in
Maronna~\cite{maronna1976}
or equations~(4.11)-(4.12) in Huber~\cite{huber1981} for M-estimators of multivariate location and covariance.
In view of this, solutions~$(\bbeta_n,\btheta_n)$ of
estimating equations~\eqref{eq:score equations M} are called M-estimators for~$(\bbeta_0,\btheta_0)$.
The expressions for $\sigma_1$ and $\sigma_2$ in Theorem~\ref{th:expansion} then
coincide with the ones in Example~3 in Tyler~\cite{tyler1982}.

As a special case, this includes the estimating equations that correspond to
maximum likelihood estimators based on independent observations
$(\by_1,\bX_1),\ldots,(\by_n,\bX_n)$ from an elliptical density~\eqref{eq:elliptical}.
The maximum likelihood estimators $(\bbeta_n,\btheta_n)$ then satisfy
estimating equations~\eqref{eq:score equations M}, for $w_1(s)=w_2(s)=-2g'(s^2)/g(s^2)$
and $w_3(s)=1$.
The expressions for $\sigma_1$ and $\sigma_2$ in Theorem~\ref{th:expansion} then
coincide with the ones in Example~2 in Tyler~\cite{tyler1982}.
\end{example}

\begin{example}[S-estimators]
\label{example:S-estimators}
S-estimators for $(\bbeta_0,\btheta_0)$ are defined  by means of a function $\rho:\R\to[0,\infty)$,
as the solution to minimizing $|\bV(\btheta)|$,
subject to
\[
\frac1n
\sum_{i=1}^n
\rho
\left(\sqrt{(\by_i-\bX_i\bbeta)^T
\bV(\btheta)^{-1}
(\by_i-\bX_i\bbeta)}\right)
\leq
b_0,
\]
where the minimum is taken over all $\bbeta\in\R^q$ and $\btheta\in\R^\ell$,
such that $\bV(\btheta)\in\text{\rm PDS}(k)$.
These estimators have been studied for linear mixed effects models in
Copt and Victoria-Feser~\cite{copt2006high}, Chervoneva and Vishnyakov~\cite{chervoneva2011,chervoneva2014}
and for general linear models with a structured covariance in
Lopuha\"a \emph{et al}~\cite{lopuhaa2023}.
According to Section~7.2 in~\cite{lopuhaa2023}, S-estimators~$(\bbeta_n,\btheta_n)$ satisfy
estimating equations~\eqref{eq:score equations M}, with $w_1(d)=\rho'(d)/d$,
$w_2(s)=k\rho'(s)/s$ and $w_3(s)=\rho'(s)s-\rho(s)+b_0$.
The expressions for $\sigma_1$ and $\sigma_2$ in Theorem~\ref{th:expansion}
coincide with the ones in Corollary~9.2 in
Lopuha\"a \emph{et al}~\cite{lopuhaa2023}.
\end{example}

\section{Homogeneous mappings of order zero}
\label{sec:homogeneous mapping}
Let $H(\bv)$ be a mapping from $\R^l$ to $\R^m$ that
is \textit{homogeneous of order zero}, that is
\begin{equation}
\label{def:homogeneous mapping}
H(\bv)=H(\alpha\bv),
\text{ for all }\alpha>0.
\end{equation}
These mappings have several applications to affine equivariant covariance estimators
that have limiting variance~\eqref{eq:limiting variance unstructured}.
Tyler~\cite{tyler1983} uses such a mapping to show that the likelihood ratio criterion
is asymptotically robust over the class of elliptical distributions.
Kent and Tyler~\cite{kent&tyler1996} consider the shape component of covariance CM-estimators
and show that the limiting variance of CM-estimators of shape depends on $\sigma_1$ only,
which may then serve as an index for the asymptotic relative efficiency.
Salibi\'an-Barrera \textit{et al}~\cite{SalibianBarrera-VanAelst-Willems2006}
derive the influence function of the shape component of covariance MM-functionals and use this
to obtain that the limiting variance of MM-estimators of shape only depends on a single scalar.
This property of the shape component is a special case of a general result in Tyler~\cite{tyler1983}
for multivariate functionals of affine equivariant covariance estimators that are asymptotically normal with
limiting variance~\eqref{eq:limiting variance unstructured}.

Estimators for a structured covariance matrix are typically not affine equivariant
and have limiting variance~\eqref{eq:limiting variance structured} instead of~\eqref{eq:limiting variance unstructured},
so that the previous results do not directly apply.
The objective of this section is to extend Theorem~1 in Tyler~\cite{tyler1983}
to estimators for a linearly structured covariance,
and discuss its consequences for corresponding estimators of shape and scale.
Moreover, we establish a similar result for estimators of the vector of variance components
and apply this to its normalized version.
We then have the following theorem.
\begin{theorem}
\label{th:limit H}
Consider $\bSigma=\bV(\btheta_0)\in\text{\rm PDS}(k)$,
for some vector $\btheta_0\in\R^\ell$
and linear variance structure $\bV$.
Let $\{\bV_n:n\geq1\}$ be a sequence of estimators for $\bSigma$
and let $\{\btheta_n:n\geq1\}$ be a sequence of estimators for the vector $\btheta_0\in\R^\ell$ of variance components.
\begin{itemize}
\item[(i)]
For $\bV\in\text{\rm PDS}(k)$, let $H(\bV)$ be continuously differentiable satisfying~\eqref{def:homogeneous mapping}.
When $\sqrt{n}(\bV_n-\bSigma)$ converges in distribution to a random matrix $\bM$,
that has a multivariate normal distribution with mean zero and variance given by~\eqref{eq:limiting variance structured},
then $\sqrt{n}(H(\bV_n)-H(\bSigma))$ is asymptotically normal with mean zero
and variance
\[
2\sigma_1 H'(\bSigma)\bL\Big(\bL^T\left(\bSigma^{-1}\otimes\bSigma^{-1}\right)\bL\Big)^{-1}\bL^TH'(\bSigma)^T.
\]
\item[(ii)]
When $\sqrt{n}(\btheta_n-\btheta_0)$ is asymptotically normal with mean zero and variance~\eqref{eq:limiting variance theta structured}.
Then for any mapping $H(\btheta)$ that satisfies~\eqref{def:homogeneous mapping},
it holds that $\sqrt{n}(H(\btheta_n)-H(\btheta_0))$ is asymptotically normal with mean zero
and variance
\[
2\sigma_1 H'(\btheta_0)
\Big(\bL^T\left(\bSigma^{-1}\otimes\bSigma^{-1}\right)\bL\Big)^{-1}
H'(\btheta_0)^T.
\]
\end{itemize}
\end{theorem}
When $\bSigma=\bV(\btheta_0)$ is unstructured, then $\vc(\bSigma)=\bL\btheta_0$,
with $\btheta_0=\text{\rm vech}(\bSigma)$, as defined in~\eqref{def:vech}, and
$\bL$ is the duplication matrix~$\mathcal{D}_k$.
Because $\bK_{k,k}H'(\bV)^T=H'(\bV)^T$, for symmetric $\bV$,
from~\eqref{eq:property duplication matrix} it follows that Theorem~\ref{th:limit H}(i)
with $\bL=\mathcal{D}_k$ recovers Theorem~1 in Tyler~\cite{tyler1983}.

From Theorem~\ref{th:limit H} it follows immediately that
the asymptotic relative efficiency of different estimators $H(\bV_n)$ for $H(\bSigma)$
can be compared by simply comparing the values of the corresponding scalar $\sigma_1$.
Similarly, the scalar $\sigma_1$ can also be used as an index for the asymptotic relative efficiency
of different estimators $H(\btheta_n)$ for $H(\btheta_0)$.
We discuss some examples below.

\begin{example}[Shape and scale of a structured covariance]
\label{ex:shape}
Suppose that $\sqrt{n}(\bV_n-\bSigma)$ is asymptotically normal
with mean zero and variance given by~\eqref{eq:limiting variance structured}.
Consider the shape component $H(\bC)=\vc(\bC)/|\bC|^{1/k}$,
where $\bC\in\text{\rm PDS}(k)$.
We have that
\begin{equation}
\label{eq:derivative shape H}
H'(\bC)
=
\frac{\partial H(\bC)}{\partial \vc(\bC)^T}
=
-
\frac1k|\bC|^{-1/k}
\text{\rm vec}(\bC)\text{\rm vec}(\bC^{-1})^T
+
|\bC|^{-1/k}
\bI_{k^2}.
\end{equation}
Then, according to Theorem~\ref{th:limit H}(i), for the shape component it follows that
$\sqrt{n}(H(\bV_n)-H(\bSigma))$ is asymptotically normal with mean zero
and variance (see Appendix for details)
\begin{equation}
\label{eq:asymp variance shape}
\frac{2\sigma_1}{|\bSigma|^{2/k}}
\left\{
\bL\Big(\bL^T\left(\bSigma^{-1}\otimes\bSigma^{-1}\right)\bL\Big)^{-1}\bL^T
-
\frac1k
\text{\rm vec}(\bSigma)
\text{\rm vec}(\bSigma)^T
\right\}.
\end{equation}
When $\bSigma$ is unstructured, then $\vc(\bSigma)=\bL\btheta_0$,
with $\btheta_0=\text{\rm vech}(\bSigma)$ and
$\bL$ is the duplication matrix~$\mathcal{D}_k$.
In that case, from~\eqref{eq:property duplication matrix} it follows that~\eqref{eq:asymp variance shape}
with $\bL=\mathcal{D}_k$
reduces to
\[
\frac{\sigma_1}{|\bSigma|^{2/k}}
\left\{
\left(
\bI_{k^2}+\bK_{k,k}
\right)(\bSigma\otimes\bSigma)
-
\frac2k\text{\rm vec}(\bSigma)\text{\rm vec}(\bSigma)^T
\right\}.
\]
This coincides with expression~(9) found in~\cite{SalibianBarrera-VanAelst-Willems2006}.
For completeness, consider the scale component $\sigma(\bC)=|\bC|^{1/(2k)}$.
It can be seen that
\begin{equation}
\label{eq:derivative scale}
\sigma'(\bC)
=
\frac1{2k}|\bC|^{1/(2k)}
\text{\rm vec}(\bC^{-1})^T.
\end{equation}
Application of the delta method then yields that $\sqrt{n}(\sigma(\bV_n)-\sigma(\bSigma))$
is asymptotically normal with mean zero and variance
\[
\frac14
\left(
\frac{2\sigma_1}{k}+\sigma_2
\right)
|\bSigma|^{1/k}.
\]
\end{example}

\begin{example}[Direction of the vector of variance components]
\label{ex:direction}
In order to create a single scalar as an index of the asymptotic efficiency for estimators $\btheta_n$
for the vector $\btheta_0$ of variance components, it is helpful to separate $\btheta_0$ into its direction and length.
The direction component $H(\btheta)=\btheta/\|\btheta\|$ satisfies~\eqref{def:homogeneous mapping}.
Its derivative is given by
\begin{equation}
\label{eq:derivative direction}
H'(\btheta)
=
\frac{\partial H(\btheta)}{\partial \btheta^T}
=
\frac{1}{\|\btheta\|}
\left(
\bI_{\ell}-\frac{\btheta\btheta^T}{\|\btheta\|^2}
\right).
\end{equation}
Then, according to Theorem~\ref{th:limit H}(ii), for the direction estimator it follows that
$\sqrt{n}(H(\btheta_n)-H(\btheta))$ is asymptotically normal with mean zero
and variance
\[
\frac{2\sigma_1}{\|\btheta_0\|^2}
\left(
\bI_{\ell}-\frac{\btheta_0\btheta_0^T}{\|\btheta_0\|^2}
\right)
\Big(\bL^T\left(\bSigma^{-1}\otimes\bSigma^{-1}\right)\bL\Big)^{-1}
\left(
\bI_{\ell}-\frac{\btheta_0\btheta_0^T}{\|\btheta_0\|^2}
\right).
\]
It does not seem possible to simplify this expression any further,
but it illustrates that one can use the scalar~$\sigma_1$ as an index for the asymptotic relative efficiency
of estimators $H(\btheta_n)$
for $H(\btheta_0)$.

An alternative is the mapping $H(\btheta)=\btheta/|\bV(\btheta)|^{1/k}$.
Since $\bV$ is linear, this $H$ also satisfies~\eqref{def:homogeneous mapping}.
For $\bV_n=\bV(\btheta_n)$, it holds that $\btheta_n=(\bL^T\bL)^{-1}\bL^T\vc(\bV_n)$,
so that
\[
H(\btheta_n)=(\bL^T\bL)^{-1}\bL^T\vc\left(\bV_n/|\bV_n|^{1/k}\right).
\]
From Example~\ref{ex:shape}, it follows that $\sqrt{n}(H(\btheta_n)-H(\btheta))$ is asymptotically normal
with mean zero and variance
\[
\frac{2\sigma_1}{|\bSigma|^{2/k}}
\left\{
\Big(\bL^T\left(\bSigma^{-1}\otimes\bSigma^{-1}\right)\bL\Big)^{-1}
-
\frac1k
\btheta_0\btheta_0^T
\right\}.
\]
This component $H$ leads to a simpler expression for the limiting variance
and the scalar $\sigma_1$ can again be used as an index for the asymptotic relative efficiency of estimators~$H(\btheta_n)$
for $H(\btheta_0)$.
\end{example}

\section{Influence function of structured covariance functionals}
\label{sec:IF}
The influence function measures the local robustness of an estimator.
It describes the effect of an infinitesimal contamination at a single point on the
corresponding functional
(see Hampel~\cite{hampel1974}).
Good local robustness is therefore illustrated by a bounded influence function.
It is defined as follows.
Let $P$ be a distribution on $\R^k$.
For $0<h<1$ and $\by\in\R^k$ fixed, define the perturbed probability measure
$P_{h,\by}=(1-h)P+h\delta_{\by}$,
where $\delta_{\by}$ denotes the Dirac measure at $\by\in\R^k$.
The \emph{influence function} of a $k\times k$ covariance functional~$\bC(\cdot)$ at probability measure $P$,
is defined as
\begin{equation}
\label{def:IF}
\text{IF}(\by;\bC,P)
=
\lim_{h\downarrow0}
\frac{\bC((1-h)P+h\delta_{\by})-\bC(P)}{h},
\end{equation}
if this limit exists.

Let $P$ be a distribution on $\R^k$ with density
$|\bSigma|^{-1/2}g\left((\by-\bmu)^T\bSigma^{-1}(\by-\bmu)\right)$, where $\bmu\in\R^k$ and $\bSigma\in\text{\rm PDS}(k)$,
and let $\bC$ be \textit{Fisher consistent} for $\bSigma$, that is $\bC(P)=\bSigma$, and \textit{affine equivariant},
meaning $\bC(P_{\bA\by+\bb})=\bA\bC(P_{\by})\bA^T$,
for any nonsingular $k\times k$ matrix $\bA$ and $\bb\in\R^k$, where $P_{\by}$ denotes the distribution of a random vector $\by$.
Croux and Haesbroeck~\cite{croux-haesbroeck2000} show that the influence function of such covariance functionals
at the $N_k(\bmu,\bSigma)$ distribution is given by
\begin{equation}
\label{eq:IF for affine C}
\text{IF}(\by;\bC,P)
=
\alpha_C(d(\by))
(\by-\bmu)(\by-\bmu)^T
-
\beta_C(d(\by))\bSigma,
\end{equation}
for some real valued functions $\alpha_C$ and $\beta_C$ and where
$d(\by)^2=(\by-\bmu)^T\bSigma^{-1}(\by-\bmu)$.
For more details on $\alpha_C$ and $\beta_C$ for different covariance functionals,
see Croux and Haesbroeck~\cite{croux-haesbroeck2000}.

Structured covariance functionals $\bM(\cdot)=\bV(\btheta(\cdot))$ are not necessarily affine equivariant,
so that the above characterizations do not directly apply.
However, Lopuha\"a \emph{et al}~\cite{lopuhaa-gares-ruizgazen2023} find similar expressions
for the influence function of the covariance S-functionals $\bM(\cdot)$ and $\btheta(\cdot)$
in a linear model with a linearly structured covariance $\bV$,
see Corollary~8.4 in~\cite{lopuhaa-gares-ruizgazen2023}.
The next lemma shows that these expressions will always appear at elliptical distributions
for covariance functionals that are a projection of some affine equivariant covariance functional.
\begin{lemma}
\label{lem:char IF}
Let $P$ be a distribution on $\R^k$ with density
$|\bSigma|^{-1/2}
g\left(
(\by-\bmu)^T\bSigma^{-1}(\by-\bmu)
\right)$, where $\bmu\in\R^k$ and $\bSigma\in\text{\rm PDS}(k)$.
Let $\bC$ be an affine equivariant covariance functional
which possesses an influence function and is Fisher consistent for $\bSigma$.
Suppose that $\bSigma=\bV(\btheta_0)$, for some $\btheta_0\in\R^\ell$, and that $\bV$ is linear
such that~$\bL$, as defined in~\eqref{def:L}, is of full column rank.
Let $\Pi_L$ be the projection matrix defined in~\eqref{def:Projection}
and define the covariance functional~$\bM$ by $\vc(\bM)=\Pi_L\vc(\bC)$.
Then the following holds.
\begin{itemize}
\item[(i)]
The functional $\bM$ is Fisher consistent for $\bSigma$ and there exist functions $\alpha_C,\beta_C:[0,\infty)\to\R$, such that
$\text{\rm IF}(\by;\vc(\bM),P)$ is given by
\[
\alpha_C(d(\by))
\bL\Big(\bL^T(\bSigma^{-1}\otimes\bSigma^{-1})\bL\Big)^{-1}
\bL^T
\vc\left(\bSigma^{-1}(\by-\bmu)(\by-\bmu)^T\bSigma^{-1}\right)
-
\beta_C(d(\by))
\vc(\bSigma),
\]
where $d^2(\by)=(\by-\bmu)^T\bSigma^{-1}(\by-\bmu)$.
\item[(ii)]
If $\btheta(P)\in\R^\ell$ is the functional, such that $\vc(\bM(\cdot))=\bL\btheta(\cdot)$, then
$\btheta$ is Fisher consistent for~$\btheta_0$ and
$\text{\rm IF}(\by;\btheta,P)$ is given by
\[
\alpha_C(d(\by))
\Big(\bL^T(\bSigma^{-1}\otimes\bSigma^{-1})\bL\Big)^{-1}
\bL^T
\vc\left(\bSigma^{-1}(\by-\bmu)(\by-\bmu)^T\bSigma^{-1}\right)
-
\beta_C(d(\by))
\btheta_0.
\]
\end{itemize}
\end{lemma}
Note that the functions $\alpha_C$ and $\beta_C$ have nothing to do with the projection $\Pi_L$,
but are inherited from the influence function~\eqref{eq:IF for affine C} of the affine equivariant
covariance functional~$\bC$.
At a distribution~$P$ that has an elliptical density~\eqref{eq:elliptical} with a linearly structured covariance,
Lopuha\"a \emph{et al}~\cite{lopuhaa-gares-ruizgazen2023} find expressions
similar to the ones in Lemma~\ref{lem:char IF} for the covariance S-functionals.
If the S-functional is defined by some function $\rho$ and constant $b_0$ (see Example~\ref{example:S-estimators}),
then
\begin{equation}
\label{def: alphaC and betaC}
\begin{split}
\alpha_C(s)
&=
\frac{k\rho'(s)}{s\delta_1}\\
\beta_C(s)
&=
\frac{\rho'(s)s}{\delta_1}
-
\frac{2(\rho(s)-b_0)}{\delta_2},
\end{split}
\end{equation}
where 
\begin{equation}
\label{def:delta12}
\begin{split}
\delta_1
&=
\frac{\mathbb{E}_{\mathbf{0},\mathbf{I}}
\left[
\rho''(\|\bz\|)\|\bz\|^2+(k+1)\rho'(\|\bz\|)\|\bz\|
\right]}{k+2}\\
\delta_2
&=
\mathbb{E}_{\mathbf{0},\mathbf{I}}
\left[
\rho'(\|\bz\|)\|\bz\|
\right].
\end{split}
\end{equation}
These $\alpha_C$ and $\beta_C$ are the same as the ones
that appear in the expression for the influence function of the affine equivariant covariance S-functional~$\bC$ in the multivariate location-scale model,
see Lopuha\"a~\cite{lopuhaa1989} or Salibi\'an-Barrera \textit{et al}~\cite{SalibianBarrera-VanAelst-Willems2006},
or in the multivariate regression model, see Van Aelst and Willems~\cite{vanaelst&willems2005}.
Indeed, the influence function $\text{IF}(\by,\vc(\bV(\btheta)),P)$ of the structured covariance functional
in Lopuha\"a~\textit{et al}~\cite{lopuhaa-gares-ruizgazen2023}
is precisely the projection $\Pi_L$ of $\text{IF}(\by,\vc(\bC),P)$ as obtained in~\cite{lopuhaa1989,
SalibianBarrera-VanAelst-Willems2006,vanaelst&willems2005}.

When $\bSigma=\bV(\btheta_0)$ is unstructured,
then $\vc(\bSigma)=\bL\btheta_0$ with $\btheta_0=\text{\rm vech}(\bSigma)$
and $\bL$ is the duplication matrix~$\mathcal{D}_k$.
In that case, from~\eqref{eq:property duplication matrix} it follows that
the expression for $\text{\rm IF}(\by;\vc(\bM),P)$ in Lemma~\ref{lem:char IF}(i) with
$\bL=\mathcal{D}_k$ reduces to
\begin{equation}
\label{eq:croux&haesbroeck}
\text{\rm IF}(\by;\vc(\bM),P)
=
\vc\left\{
\alpha_C(d(\by))(\by-\bmu)(\by-\bmu)^T
-
\beta_C(d(\by))\bSigma
\right\}.
\end{equation}
This coincides with the expression found in Lemma~1 in
Croux and Haesbroeck~\cite{croux-haesbroeck2000}.

Mappings $H$ that satisfy~\eqref{def:homogeneous mapping} also have useful applications
to influence functions of affine equivariant covariance functionals $\bC$ and their the gross-error-sensitivity (GES).
Kent and Tyler~\cite{kent&tyler1996} consider functionals $\bC/|\bC|^{1/k}$
and $\bC/\text{tr}(\bC)$ to obtain that the GES of different CM-functionals
is proportional to a single scalar.
Salibi\'an-Barrera \textit{et al}~\cite{SalibianBarrera-VanAelst-Willems2006}
derive the influence function of the shape component of covariance MM-functionals and
show that it is proportional to a single function~$\alpha_C$,
which no longer depends on the scale-functional used in the first step.
In fact, these properties hold more general for functionals $H$ satisfying~\eqref{def:homogeneous mapping}
applied to affine equivariant covariance functionals.
The next lemma establishes similar results for linearly structured covariance functionals.
\begin{lemma}
\label{lem:IF for H}
Let $P$ be a distribution on $\R^k$ with an elliptical contoured density~\eqref{eq:elliptical}.
Suppose that $\bSigma=\bV(\btheta_0)$, for  some $\btheta_0\in\R^\ell$, and that $\bV$ is linear
such that~$\bL$, as defined in~\eqref{def:L}, is of full column rank.
\begin{itemize}
\item[(i)]
Let~$\bM\in\text{\rm PDS}(k)$ be a covariance functional
that is Fisher consistent for $\bSigma$ and
which possesses an influence function given by
Lemma~\ref{lem:char IF}(i).
Let $H(\bM)$ be continuously differentiable in a neighborhood of~$\bM(P)$
satisfying~\eqref{def:homogeneous mapping}.
Then $\text{\rm IF}(\by;H(\bM),P)$ is given by
\[
\alpha_C(d(\by))
H'(\bSigma)
\bL\Big(\bL^T(\bSigma^{-1}\otimes\bSigma^{-1})\bL\Big)^{-1}
\bL^T
\left(\bSigma^{-1}\otimes\bSigma^{-1}\right)
\Big((\by-\bmu)\otimes(\by-\bmu)\Big),
\]
where $d^2(\by)=(\by-\bmu)^T\bSigma^{-1}(\by-\bmu)$.
\item[(ii)]
Let $\btheta\in\R^{\ell}$ be a functional that is Fisher consistent for $\btheta_0$
and which possesses an influence function
given by Lemma~\ref{lem:char IF}(ii).
Let $H(\btheta)$ be continuously differentiable in a neighborhood of~$\btheta(P)$
satisfying~\eqref{def:homogeneous mapping}.
Then $\text{\rm IF}(\by;H(\btheta),P)$ is given by
\[
\alpha_C(d(\by))
H'(\btheta_0)
\Big(\bL^T(\bSigma^{-1}\otimes\bSigma^{-1})\bL\Big)^{-1}
\bL^T
\left(\bSigma^{-1}\otimes\bSigma^{-1}\right)
\Big((\by-\bmu)\otimes(\by-\bmu)\Big).
\]
\end{itemize}
\end{lemma}

Consider the GES defined by $\sup_{\by\in\R^k}\|\text{IF}(\by;\cdot)\|$, for some norm $\|\cdot\|$.
From Lemma~\ref{lem:IF for H} it follows immediately that regardless of the choice of the norm, the value $\|\text{IF}(\by;H(\bM),P)\|$
for different functionals $H(\bM(P))$
is proportional to $|\alpha_C(d(\by))|$ and similarly for functionals $H(\btheta(P))$.
We discuss some examples below.

\begin{example}[Shape and scale of a structured covariance]
\label{ex:shape for M}
For the shape functional $H(\bM)=\vc(\bM)/|\bM|^{1/k}$,
from Lemma~\ref{lem:IF for H}(i) together with~\eqref{eq:derivative shape H}
we find
\[
\text{\rm IF}(\by;H(\bM),P)
=
-
\frac1k|\bSigma|^{-1/k}
\text{\rm tr}
\left(
\bSigma^{-1}
\text{\rm IF}(\by;\bM,P)
\right)\cdot
\vc(\bSigma)
+
|\bSigma|^{-1/k}
\text{\rm IF}(\by;\vc(\bM),P).
\]
See also Salibi\'an \emph{et al}~\cite{SalibianBarrera-VanAelst-Willems2006}.
In particular, at a distribution $P$ with an elliptically contoured density with parameters~$\bmu$ and~$\bSigma=\bV(\btheta_0)$
one finds that~$\text{\rm IF}(\by;H(\bM),P)$ is given by
\begin{equation}
\label{eq:IF shape cov}
\begin{split}
\frac{\alpha_C(d(\by))}{|\bSigma|^{1/k}}
\bigg\{
\bL\Big(\bL^T(\bSigma^{-1}\otimes\bSigma^{-1})\bL\Big)^{-1}
\bL^T
\vc\left(\bSigma^{-1}(\by-\bmu)(\by-\bmu)^T\bSigma^{-1}\right)
-
\frac{d(\by)^2}k
\vc(\bSigma)
\bigg\},
\end{split}
\end{equation}
where $d(\by)^2=(\by-\bmu)^T\bSigma^{-1}(\by-\bmu)$.
It follows that $\|\text{\rm IF}(\by;H(\btheta),P)\|$
will be proportional
to $|\alpha_C(d(\by))d(\by)^2|$.
When $\bSigma$ is unstructured, then $\vc(\bSigma)=\bL\btheta_0$,
where $\btheta_0=\text{\rm vech}(\bSigma)$, as defined in~\eqref{def:vech}, and
$\bL$ is the duplication matrix~$\mathcal{D}_k$.
In that case, from~\eqref{eq:property duplication matrix} it follows that~\eqref{eq:IF shape cov}
with $\bL=\mathcal{D}_k$
reduces to
\[
\frac{\alpha_C(d(\by))}{|\bSigma|^{1/k}}
\vc
\left\{
(\by-\bmu)(\by-\bmu)^T-\frac{d(\by)^2}{k}\bSigma
\right\},
\]
which coincides with formula~(3) in~\cite{SalibianBarrera-VanAelst-Willems2006}.
For completeness, consider the scale component $\sigma(\bM)=|\bM|^{1/(2k)}$.
From~\eqref{eq:derivative scale}, it follows that
\[
\text{\rm IF}(\by;\sigma,P)
=
\frac12|\bSigma|^{-1/(2k)}
\gamma_C(d(\by)),
\]
where $\gamma_C(s)
=
\alpha_C(s)
s^2/k
-
\beta_C(s)$,
which matches with equation~(4) in~\cite{SalibianBarrera-VanAelst-Willems2006}.
\end{example}

\begin{example}[Direction of the vector of variance components]
\label{ex:direction IF}
For the direction functional $H(\btheta)=\btheta/\|\btheta\|$,
from Lemma~\ref{lem:IF for H}(ii) together with~\eqref{eq:derivative direction} we find that,
at a distribution $P$ with an elliptically contoured distribution with parameters $\bmu$ and $\bSigma=\bV(\btheta_0)$,
$\text{\rm IF}(\by;H(\btheta),P)$ is given by
\[
\alpha_C(d(\by))
\left(
\frac{1}{\|\btheta_0\|}
\bI_\ell
-
\frac{\btheta_0\btheta_0^T}{\|\btheta_0\|^3}
\right)
\Big(\bL^T(\bSigma^{-1}\otimes\bSigma^{-1})\bL\Big)^{-1}
\bL^T
\vc\left(\bSigma^{-1}(\by-\bmu)(\by-\bmu)^T\bSigma^{-1}\right).
\]
It follows that $\|\text{\rm IF}(\by;H(\btheta),P)\|$
will be proportional
to $|\alpha_C(d(\by))d(\by)^2|$.
An alternative is the mapping $H(\btheta)=\btheta/|\bV(\btheta)|^{1/k}$.
Since $\bV$ is linear, $H$ satisfies~\eqref{def:homogeneous mapping}.
For $\bM(P)=\bV(\btheta(P))$, it holds that $\btheta(P)=(\bL^T\bL)^{-1}\bL^T\vc(\bM(P))$,
so that
\[
H(\btheta(P))=(\bL^T\bL)^{-1}\bL^T\vc(\bM(P))/|\bM(P)|^{1/k}.
\]
From Example~\ref{ex:shape for M} it follows that
$\text{\rm IF}(\by;H(\btheta),P)$ is given by
\[
\begin{split}
&
\frac{\alpha_C(d(\by))}{|\bSigma|^{1/k}}
\bigg\{
\Big(\bL^T(\bSigma^{-1}\otimes\bSigma^{-1})\bL\Big)^{-1}
\bL^T
\vc\left(\bSigma^{-1}(\by-\bmu)(\by-\bmu)^T\bSigma^{-1}\right)
-
\frac{d(\by)^2}k
\btheta_0
\bigg\}
\end{split}
\]
using that $\vc(\bSigma)=\bL\btheta_0$.
Again we find that $\|\text{\rm IF}(\by;H(\btheta),P)\|$ is proportional to 
$|\alpha_C(d(\by))d(\by)^2|$.
\end{example}

\section{Application}
\label{sec:application}
We apply our results to S-estimators and S-functionals in the linear model~\eqref{def:linear model}.
Let $P$ be the distribution for the random variable $\bs=(\by,\bX)$, which is such that $\by\mid\bX$
has an elliptically contoured distribution~\eqref{eq:elliptical} with parameters
$\bmu=\bX\bbeta_0$ and $\bSigma=\bV(\btheta_0)=\theta_{01}\bL_1+\cdots+\theta_{0\ell}\bL_\ell$.
Consider the S-estimator for $(\bbeta_0,\btheta_0)$
defined as the solution to minimizing $|\bV(\btheta)|$,
subject to
\[
\frac1n
\sum_{i=1}^n
\rho
\left(\sqrt{(\by_i-\bX_i\bbeta)^T
\bV(\btheta)^{-1}
(\by_i-\bX_i\bbeta)}\right)
=
b_0,
\]
where the minimum is taken over all $\bbeta\in\R^q$ and $\btheta\in\R^\ell$,
such that $\bV(\btheta)\in\text{\rm PDS}(k)$.
For the function $\rho$ we take Tukey's bi-weight
\begin{equation}
\label{def:biweight}
\rho_{\mathrm{B}}(s;c)
=
\begin{cases}
s^2/2-s^4/(2c^2)+s^6/(6c^4), & |s|\leq c;\\
c^2/6 & |s|>c, 
\end{cases}
\end{equation}
and $b_0=\E_{\mathbf{0},\bI_k}[\rho_B(\|\bz\|;c)]$.
From Theorem~6.1 in Lopuha\"a~\textit{\textit{et al}}~\cite{lopuhaa2023}
it is known that the breakdown point of the S-estimator depends on the cut-off constant $c$
and is at least $\lceil nb_0/(c^2/6)\rceil/n$, or asymptotically $\epsilon^*=b_0/(c^2/6)$.
\begin{table}
\caption{Cut-off values of $\rho_B$ for different breakdown points and dimensions.}
\label{tab:cutoff}
\centering
\[
\begin{array}{crrrrrrrrrr}
\hline
\hline
&&&&&&&&&&\\[-5pt]
  & \multicolumn{10}{c}{\text{Breakdown point}} \\
&&&&&&&&&&\\[-5pt]
k & 0.05    & 0.10   &  0.15  & 0.20   & 0.25   & 0.30   & 0.35    & 0.40  &  0.45  & 0.50\\
   \cline{2-11}
&&&&&&&&&&\\[-5pt]
1 & 7.545   &  5.182 &  4.096 &  3.421 &  2.937 &  2.561 &  2.252 &  1.988 &  1.756 &  1.548\\
2 &  10.767 &  7.474 &  5.981 &  5.069 &  4.427 &  3.938 &  3.542 &  3.209 &  2.920 &  2.661\\
5 &  17.114 & 11.950 &  9.628 &  8.220 &  7.242 &  6.505 &  5.918 &  5.432 &  5.017 &  4.652\\
10 & 24.246 & 16.961 & 13.694 & 11.719 & 10.351 &  9.324 &  8.510 &  7.840 &  7.271 &  6.776\\
&&&&&&&&&& \\[-5pt]
\hline
\hline
\end{array}
\]
\end{table}
Table~\ref{tab:cutoff} gives the cut-off values of $\rho_B$ for given asymptotic lower bounds
$\epsilon^*=0.05,0.10,\ldots,0.50$ on the breakdown point
in dimensions $k=1,2,5,10$.
This table partly overlaps with Table~3 in Rousseeuw and Yohai~\cite{rousseeuw-yohai1984}.

According to Corollary~9.2 in Lopuha\"a~\textit{et al}~\cite{lopuhaa-gares-ruizgazen2023},
the scalar
$\lambda=\E_{\mathbf{0},\bI_k}\left[\rho_B'(\|\bz\|;c)^2\right]/(k\alpha^2)$
represents the asymptotic efficiency of the regression S-estimator $\bbeta_n$ relative to the least squares estimator
(for which $\lambda=1)$,
where
\begin{equation}
\label{def:alpha}
\alpha
=
\mathbb{E}_{\mathbf{0},\bI_k}
\left[
\left(1-\frac{1}{k}\right)
\frac{\rho_B'(\|\bz\|;c)}{\|\bz\|}
+
\frac1k
\rho_B''(\|\bz\|;c)
\right].
\end{equation}
From Examples~\ref{ex:shape} and~\ref{ex:direction}, together with Theorem~\ref{th:expansion} and Example~\ref{example:S-estimators},
it follows that the scalar
\[
\sigma_1
=
\frac{k\E_{\mathbf{0},\bI_k}\left[\rho_B'(\|\bz\|;c)^2\|\bz\|^2\right]}{(k+2)\delta_1^2},
\]
where $\delta_1$ is defined in~\eqref{def:delta12},
serves as an index for the asymptotic efficiency of both the S-estimator 
of shape as well as the S-estimator for the direction of the vector of variance components,
relative to the least squares estimators of shape and direction, respectively (for which $\sigma_1=1$).
Finally, from Example~\ref{ex:shape}, together with Example~\ref{example:S-estimators},
it follows that
\[
\sigma_3
=
\frac14
\left(
\frac{2\sigma_1}{k}+\sigma_2
\right)
=
\frac{\mathbb{E}_{\textbf{0},\mathbf{I}_k}
\left[
\left(
\rho_B(\|\bz\|;c)-b_0
\right)^2
\right]}{\delta_2^2},
\]
where $\delta_2$ is defined in~\eqref{def:delta12},
serves as an index for the asymptotic efficiency of the S-estimator of scale
relative the least squares (for which $\sigma_3=1/(2k)$).
As a consequence, the cutoff constant~$c$ of $\rho_B$ can be tuned in such a way that the asymptotic efficiency
$1/\lambda$ relative to the least squares estimator is high at the normal distribution and similarly for $1/\sigma_1$
and $1/(2k\sigma_3)$.
Since $c$ also determines the breakdown point,
this forces a trade-off between efficiency and breakdown point.
Typically, large values of $c$ correspond to high efficiency and low breakdown point,
and vice-versa for moderate values of $c$.

We further investigate how this trade-off relates to the gross error sensitivity (GES) of the corresponding S-functionals.
For simplicity we only consider perturbations in $\by$ and leave $\bX$ unchanged.
From Corollary~8.4 in Lopuha\"a~\textit{et al}~\cite{lopuhaa-gares-ruizgazen2023},
for the regression S-functional it then follows that $\|\text{IF}(\by;\bbeta,P)\|$
is proportional to $\alpha^{-1}\left|\rho_B'(d(\by);c)\right|$,
where $\alpha$ is defined in~\eqref{def:alpha} and $d(\by)^2=(\by-\bX\bbeta_0)^T\bSigma^{-1}(\by-\bX\bbeta_0)$.
Therefore, we propose the scalar
\[
G_1
=
\frac{1}{\alpha}
\sup_{s>0}
\left|
\rho_B'(s;c)
\right|,
\]
as an index for the GES of regression S-functionals.
This coincides with the GES index for location CM-functionals in Kent and Tyler~\cite{kent&tyler1996}.
From Examples~\ref{ex:shape for M} and~\ref{ex:direction IF}, together with Lemma~\ref{lem:IF for H} and~\eqref{def: alphaC and betaC},
for both the shape and direction S-functional, it follows that~$\|\text{IF}(\by)\|$
is proportional to $\delta_1^{-1}|\rho_B'(d(\by);c)d(\by)|$, where
$\delta_1$ is defined in~\eqref{def:delta12}.
We propose
the scalar
\[
G_2
=
\frac{k}{(k+2)\delta_1}
\sup_{s>0}
\left|
\rho_B'(s;c)s
\right|,
\]
as an index for the GES of shape and direction S-functionals.
In this way, $G_2$ coincides with the GES index for CM-functionals of shape in Kent and Tyler~\cite{kent&tyler1996}.
Finally, from Example~\ref{ex:shape for M} and~\eqref{def: alphaC and betaC}, if follows that 
for the scale functional $\|\text{IF}(\by)\|$ is proportional to 
$\delta_2^{-1}|\rho_B(d(\by);c)-b_0|$, where
$\delta_2$ is defined in~\eqref{def:delta12}.
We propose 
\[
G_3=\frac{1}{\delta_2}
\sup_{s>0}
\left|
\rho_B(s;c)-b_0
\right|,
\]
as an index for the GES of the S-functional of scale.

\begin{figure}[t]
  \centering
  \includegraphics[width=\textwidth]{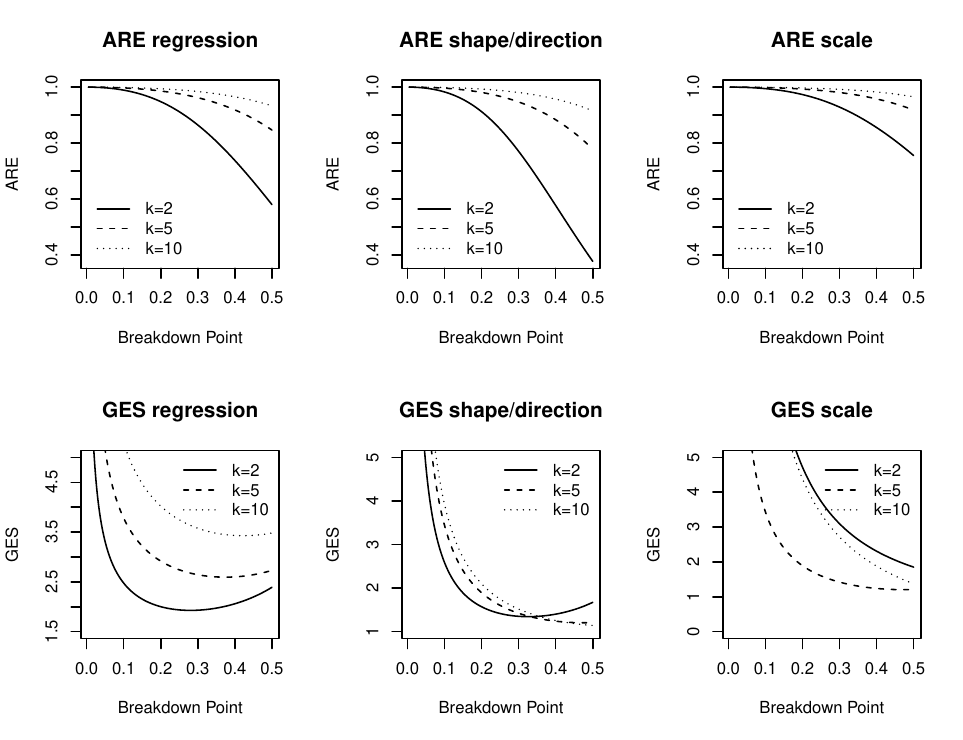}
  \caption{ARE and GES as functions of the breakdown point at the multivariate normal in dimensions~$k=2,5,10$.}
  \label{fig:dim2}
\end{figure}

We investigate how the asymptotic efficiency at the normal distribution of the S-estimators,
and the GES of the corresponding S-functionals behave as we vary the breakdown point of the 
S-estimator between 0 and 0.5.
Given a value $\epsilon^*$ of the breakdown point, we determine the corresponding cut-off constant $c$
by solving $\epsilon^*=\E_{\mathbf{0},\bI_k}[\rho_B(\|\bz\|;c)]/(c^2/6)$.
With this value of $c$, we compute the values of $\lambda$, $\sigma_1$ and $\sigma_3$
and the GES indices $G_1$, $G_2$ and $G_3$.
In Figure~\ref{fig:dim2}, on the top row we have plotted the asymptotic relative efficiencies $1/\lambda$, $1/\sigma_1$ and $1/(2k\sigma_3)$
as a function of the breakdown point for dimensions $k=2,5,10$,
and the bottom row contains plots of the GES indices $G_1$, $G_2$ and $G_3$ for the same dimensions.
As expected, the efficiency decreases with increasing breakdown point, but
the loss of efficiency is less severe for the S-estimator of scale compared to 
the S-estimator for regression and the S-estimators for shape and direction.
In dimension $k=2$ (solid lines), the 50\% breakdown S-estimators have asymptotic efficiencies
$1/\lambda=0.580$, $1/\sigma_1=0.376$, and $1/(4\sigma_3)=0.755$.
However, one can gain both efficiency and lower the GES at the cost of a lower breakdown point.
For example, the GES index of the regression functional attains its minimal value $G_1=1.927$ at breakdown point 28\%,
which corresponds to cut-off value $c=4.115$.
For this cut-off value the GES index of the shape and direction functional is $G_2=1.368$, which is not far off from its minimal value 1.344,
and the GES index for scale is $G_3=3.323$.
Furthermore, the asymptotic efficiencies then become $1/\lambda=0.884$, $1/\sigma_1=0.803$,
and $1/(4\sigma_3)=0.939$, for the regression estimator, the estimators of shape and direction, and the scale estimator,
respectively.
Similarly, the GES index of the shape and direction functionals attains its minimal value $G_2=1.344$ for $c=3.722$.
This would yield $G_1=1.947$, $G_3=2.844$, $1/\lambda=0.835$, $1/\sigma_1=0.723$, $1/(4\sigma_3)=0.912$ and breakdown point 33\%.
The GES index of the scale functional attains its minimum value $G_3=1.852$ at 50\% breakdown point,
so no simultaneous gain in efficiency and smaller GES values $G_1$ and $G_2$ can be achieved at the cost of a smaller breakdown point.

In dimension $k=5$ (dashed lines), the 50\% breakdown S-estimators have asymptotic efficiencies
$1/\lambda=0.864$, $1/\sigma_1=0.778$, and $1/(4\sigma_3)=0.918$.
The GES index of the regression functional attains its minimal value $G_1=2.595$ at breakdown point 37\%.
The corresponding GES index for shape and direction functionals is $G_2=1.271$ and $G_3=1.480$ for the scale functionals.
Corresponding to this smaller regression GES index we observe a gain in the asymptotic efficiencies: 
$1/\lambda=0.932$, $1/\sigma_1=0.903$,
and $1/(4\sigma_3)=0.965$, for the regression estimator, the estimators of shape and direction, and the scale estimator,
respectively.
The GES index of the shape and direction functionals attains its minimal value at breakdown point 47\%, so
the gain in both efficiency and a smaller $G_2$ value is negligible.
The situation for the GES index for scale is the same as in dimension $k=2$,
where no simultaneous gain in efficiency and smaller GES values $G_1$ and $G_2$ can be achieved at the cost of a smaller breakdown point.

Finally, in dimension (dotted lines), the 50\% breakdown S-estimators have asymptotic efficiencies
$1/\lambda=0.933$, $1/\sigma_1=0.915$, and $1/(4\sigma_3)=0.965$.
The GES index of the regression functional attains its minimal value $G_1=3.426$ at breakdown point 42\%.
The corresponding GES index for shape and direction functionals is $G_2=1.221$ and $G_3=1.744$ for the scale functionals.
Corresponding to this smaller regression GES index we observe a gain in the asymptotic efficiencies: 
$1/\lambda=0.960$, $1/\sigma_1=0.949$,
and $1/(4\sigma_3)=0.979$, for the regression estimator, the estimators of shape and direction, and the scale estimator,
respectively.
Both GES indices $G_2$ and $G_3$ attain their minimal values at 50\% breakdown,
so no simultaneous gain in efficiency and smaller GES value $G_1$ can be achieved at the cost of a smaller breakdown point.

We conclude that at a moderate loss of breakdown point, from 50\% to about 30\%-40\%, one can gain efficiency 
of the S-estimators and at the same time reduce the GES of the regression S-estimator.
The improvements becomes less as the dimension increases.

\appendix

\section{Proofs}
\paragraph{Proof of Theorem~\ref{th:projection}.}
\begin{proof}
It can be seen that the projection matrix, as defined in~\eqref{def:Projection}, is given by
\begin{equation}
\label{def:Projection matrix}
\Pi_L
=
\bL\left(\bL^T\left(\bSigma^{-1}\otimes\bSigma^{-1}\right)\bL\right)^{-1}\bL^T\left(\bSigma^{-1}\otimes\bSigma^{-1}\right).
\end{equation}
Since $\bN$ is of radial type with respect to $\bSigma$,
it follows from Corollary~1 in  Tyler~\cite{tyler1982} that there
exist constants~$\eta$, $\sigma_1$ and $\sigma_2$ with $\sigma_1\geq 0$ and $\sigma_2\geq -2\sigma_1/k$,
such that
$\E[\bN]=\eta\bSigma$ and
\[
\text{var}\{\vc(\bN)\}
=
\sigma_1(\bI_{k^2}+\bK_{k,k})(\bSigma\otimes\bSigma)+\sigma_2\vc(\bSigma)\vc(\bSigma)^T.
\]
Since $\bV$ is linear, it holds that $\vc(\bSigma)=\bL\btheta$, so that $\Pi_L\vc(\bSigma)=\vc(\bSigma)$.
It follows that $\bM$ has expectation
\[
\E[\vc(\bM)]
=
\Pi_L\vc(\E[\bN])
=
\eta
\Pi_L\vc(\bSigma)
=
\eta
\vc(\bSigma),
\]
and variance
\[
\begin{split}
\text{var}(\vc(\bM))
&=
\Pi_L
\text{var}(\vc(\bN))
\Pi_L^T\\
&=
\sigma_1
\Pi_L(\bI_{k^2}+\bK_{k,k})(\bSigma\otimes\bSigma)\Pi_L^T
+
\sigma_2
\Pi_L\vc(\bSigma)\vc(\bSigma)^T\Pi_L^T\\
&=
\sigma_1
\Pi_L(\bI_{k^2}+\bK_{k,k})(\bSigma\otimes\bSigma)\Pi_L^T
+
\sigma_2
\vc(\bSigma)\vc(\bSigma)^T.
\end{split}
\]
Note that
\begin{equation}
\label{eq:property K}
\begin{split}
\bK_{k,k}(\bA\otimes\bB)
&=
(\bB\otimes\bA)\bK_{k,k}\\
\bK_{k,k}\vc(\bA)
&=
\vc(\bA^T),
\end{split}
\end{equation}
e.g., see~\cite[Chapter~3, Section 7]{magnus&neudecker1988}.
Since $\bSigma=\bV(\btheta)$ is symmetric, also $\bL_j=\partial\bV/\partial\theta_j$ is symmetric,
for $j=1,\ldots,\ell$.
This means that $\bK_{k,k}\bL=\bL$ and it follows that
\[
\begin{split}
\Pi_L(\bI_{k^2}+\bK_{k,k})(\bSigma\otimes\bSigma)\Pi_L^T
&=
\bL\left(\bL^T\left(\bSigma^{-1}\otimes\bSigma^{-1}\right)\bL\right)^{-1}\bL^T
(\bI_{k^2}+\bK_{k,k})\Pi_L^T\\
&=
2\bL\left(\bL^T\left(\bSigma^{-1}\otimes\bSigma^{-1}\right)\bL\right)^{-1}\bL^T
\Pi_L^T\\
&=
2\bL\left(\bL^T\left(\bSigma^{-1}\otimes\bSigma^{-1}\right)\bL\right)^{-1}\bL^T.
\end{split}
\]
This finishes the proof of part~(i).
Since $\bL$ has full rank, it holds that
\begin{equation}
\label{eq:from vec to theta}
(\bL^T\bL)^{-1}\bL^T\vc(\bSigma)
=
(\bL^T\bL)^{-1}\bL^T\bL\btheta
=
\btheta,
\end{equation}
and similarly $(\bL^T\bL)^{-1}\bL^T\vc(\bM)=\bT$.
This immediately gives
\[
\E[\bT]
=
(\bL^T\bL)^{-1}\bL^T\E[\vc(\bM)]
=
\eta(\bL^T\bL)^{-1}\bL^T\vc(\bSigma)
=
\eta\btheta,
\]
and
\[
\text{\rm var}(\bT)
=
(\bL^T\bL)^{-1}\bL^T
\text{\rm var}\{\vc(\bM)\}
\bL(\bL^T\bL)^{-1}.
\]
When we insert~\eqref{eq:limiting variance structured} and apply~\eqref{eq:from vec to theta},
the theorem follows.
\end{proof}

\paragraph{Proof of Theorem~\ref{th:expansion}.}
The proof follows the line of reasoning used in the proofs of Theorem~9.1 and
Corollary~9.2 in Lopuha\"a \emph{et al}~\cite{lopuhaa-gares-ruizgazen2023} for S-estimators.
These proofs are based on estimating equations~\eqref{eq:score equations M} with
$w_1(d)=\rho'(d)/d$, $w_2(d)=k\rho'(d)/d$ and $w_3(d)=\rho'(d)d-\rho(d)+b_0$,
and require conditions (R1)-(R5) in~\cite{lopuhaa-gares-ruizgazen2023}
on the function $\rho$.
For the proof of Theorem~\ref{th:expansion} these conditions have been reformulated
into similar conditions (C1)-(C3) for general $w_1$, $w_2$, and $w_3$.
Furthermore, in order to incorporate the case $w_1=w_2=w_3=1$ of Example~\ref{example:LS},
we have slightly adapted some of the boundedness conditions and use that
\[
d^2=(\by-\bX\bbeta)^T\bV^{-1}(\by-\bX\bbeta)
\leq
\frac{\|\by-\bX\bbeta\|^2}{\lambda_k(\bV)}
\leq
\frac{(\|\by\|+\|\bX\|\cdot\|\bbeta\|)^2}{\lambda_k(\bV)}
\leq
\frac{\|\bs\|^2(1+\|\bbeta\|)^2}{\lambda_k(\bV)}.
\]
This will ensure that $d^2$ is bounded by a multiple of $\|\bs\|^2$ on a neighborhood of $\bxi_0$.
In order to apply dominated convergence, we then require $\E\|\bs\|^4<\infty$
in Theorem~\ref{th:expansion} instead of $\E\|\bX\|^2<\infty$, which was sufficient for
Corollary~9.2 in~\cite{lopuhaa-gares-ruizgazen2023}.

\begin{proof}
Define
\begin{equation}
\label{def:Lambda}
\Lambda(\bxi)
=
\int \Psi(\bs,\bxi)\,\dd P(\bs).
\end{equation}
From the properties of elliptically contoured densities, one has that $\E\left[\Psi(\bs,\bxi_0)\big|\bX\right]=\mathbf{0}$,
so that~$\Lambda(\bxi_0)=\mathbf{0}$.
Conditions (C1)-(C3) yield that $\Lambda$ is continuously differentiable
and by application of empirical process theory (see e.g., Lemma~11.8 in~\cite{supplement} for the special case of S-estimators)
one finds
\begin{equation}
\label{eq:decomposition estimator}
\begin{split}
\mathbf{0}
&=
\int \Psi(\mathbf{s},\bxi_n)\,\dd P(\mathbf{s})
+
\int \Psi(\mathbf{s},\bxi_0)\,\dd (\mathbb{P}_n-P)(\mathbf{s})
+
o_P(n^{-1/2})\\
&=
\Lambda'(\bxi_0)(\bxi_n-\bxi_0)
+
\frac{1}{n}
\sum_{i=1}^{n}
\Big\{
\Psi(\bs_i,\bxi_0)
-
\E[\Psi(\bs_i,\bxi_0)]
\Big\}
+
o_P(n^{-1/2}).
\end{split}
\end{equation}
Similar to Lemma~8.3 in Lopuha\"a \emph{et al}~\cite{lopuhaa-gares-ruizgazen2023},
we find that $\Lambda'(\bxi_0)$ is a block matrix.
This implies that $\sqrt{n}(\bbeta_n-\bbeta_0)$ and $\sqrt{n}(\btheta_n-\btheta_0)$ are
asymptotically independent
and from~\eqref{eq:decomposition estimator} we obtain
\[
\begin{split}
\sqrt{n}(\vc(\bV_n)-\vc(\bSigma))
&=
\bL\sqrt{n}(\btheta_n-\btheta_0)\\
&=
-
\bL
\Lambda_{\btheta}'(\bxi_0)^{-1}
\frac{1}{\sqrt{n}}
\sum_{i=1}^{n}
\Big\{
\Psi_{\btheta}(\bs_i,\bxi_0)
-
\E[\Psi_{\btheta}(\bs_i,\bxi_0)]
\Big\}
+
o_P(1),
\end{split}
\]
where $\Psi_{\btheta}$ is defined in~\eqref{eq:Psi function}, and
\[
\Lambda_{\btheta}'(\bxi_0)
=
\int
\frac{\partial\Psi_{\btheta}(\bs,\bxi_0)}{\partial \btheta}\,\dd P(\bs)
=
\gamma_1
\bL^T\left(\bSigma^{-1}\otimes\bSigma^{-1}\right)\bL
-
\gamma_2\bL^T
\vc(\bSigma^{-1})
\vc(\bSigma^{-1})^T
\bL,
\]
where
\begin{equation}
\label{def:gamma12}
\begin{split}
\gamma_1
&=
\frac{\E_{\mathbf{0},\bI_k}
\left[
w_2'(\|\bz\|)\|\bz\|^3+k(k+2)w_3(\|\bz\|)
\right]}{k(k+2)}\\
\gamma_2
&=
\frac{\E_{\mathbf{0},\bI_k}
\left[
(k+2)w_3'(\|\bz\|)\|\bz\|
-
w_2'(\|\bz\|)\|\bz\|^3
\right]}{2k(k+2)}.
\end{split}
\end{equation}
First note that we can write
\begin{equation}
\label{eq:decomp psi-theta}
\Psi_{\btheta}(\bs,\bxi_0)
=
\bL^T(\bSigma^{-1}\otimes\bSigma^{-1})
\vc\left\{\Psi_{\bC}(\bs,\bxi_0)\right\}
\end{equation}
where
\begin{equation}
\label{def:PsiC}
\Psi_{\bC}(\bs,\bxi_0)
=
w_2(d_0)(\by-\bX\bbeta_0)(\by-\bX\bbeta_0)^T-w_3(d_0)\bSigma,
\end{equation}
with $d_0^2=(\by-\bX\bbeta_0)^T\bSigma^{-1}(\by-\bX\bbeta_0)$.
Furthermore, similar to Lemma~11.5 in~\cite{supplement}, we find that
\[
\Lambda_{\btheta}'(\bxi_0)^{-1}
=
a
(\bE^T\bE)^{-1}
+b
(\bE^T\bE)^{-1}
\bE^T
\vc(\bI_k)
\vc(\bI_k)^T
\bE(\bE^T\bE)^{-1},
\]
where $\bE=\left(\bSigma^{-1/2}\otimes\bSigma^{-1/2}\right)\bL$
and $a=1/\gamma_1$ and $b=\gamma_2/(\gamma_1(\gamma_1-k\gamma_2))$,
with
\[
\gamma_1-k\gamma_2
=
\frac{1}{2k}
\E_{\mathbf{0},\bI_k}\Big[w_2'(\|\bz\|)\|\bz\|^3+2kw_3(\|\bz\|)-kw_3'(\|\bz\|)\|\bz\|\Big]
\ne 0.
\]
This means that $\bE^T\bE=\bL^T(\bSigma^{-1}\otimes\bSigma^{-1})\bL$
and since $\vc(\bSigma)=\bL\btheta_0$, we have
\[
(\bE^T\bE)^{-1}
\bE^T
\vc(\bI_k)
=
(\bL^T(\bSigma^{-1}\otimes\bSigma^{-1})\bL)^{-1}
\bL(\bSigma^{-1}\otimes\bSigma^{-1})
\vc(\bSigma)
=
\btheta_0,
\]
so that $\Lambda_{\btheta}'(\bxi_0)^{-1}=a(\bL^T(\bSigma^{-1}\otimes\bSigma^{-1})\bL)^{-1}+b\btheta_0\btheta_0^T$.
Furthermore, since $\vc(\bSigma)=\bL\btheta_0$, together with 
\begin{equation}
\label{eq:kronecker}
\vc(\bA\bB\bC)=(\bC^T\otimes\bA)\vc(\bB)
\end{equation}
e.g., see~\cite[Chapter~2, Section 4]{magnus&neudecker1988},
and $\Pi_L$ as in~\eqref{def:Projection matrix}, we find
\[
\begin{split}
\bL\Lambda_{\btheta}'(\bxi_0)^{-1}
\bL^T(\bSigma^{-1}\otimes\bSigma^{-1})
&=
a
\Pi_L
+
b\bL\btheta_0\btheta_0^T\bL^T(\bSigma^{-1}\otimes\bSigma^{-1})
=
a
\Pi_L
+
b\vc(\bSigma)\vc(\bSigma^{-1})^T.
\end{split}
\]
It follows that
\[
\begin{split}
\bL\Lambda_{\btheta}'(\bxi_0)^{-1}
\Psi_{\btheta}(\bs,\bxi_0)
&=
\bL\Lambda_{\btheta}'(\bxi_0)^{-1}
\bL^T(\bSigma^{-1}\otimes\bSigma^{-1})
\vc\left\{\Psi_{\bC}(\bs,\bxi_0)\right\}\\
&=
a\Pi_L\vc\left\{\Psi_{\bC}(\bs,\bxi_0)\right\}
+
b
\vc(\bSigma)
\vc(\bSigma^{-1})^T
\vc\left\{\Psi_{\bC}(\bs,\bxi_0)\right\}\\
&=
a\Pi_L\vc\left\{\Psi_{\bC}(\bs,\bxi_0)\right\}
+
b
\vc(\bSigma)
\text{tr}
\left\{
\bSigma^{-1}\Psi_{\bC}(\bs,\bxi_0)
\right\}\\
&=
a\Pi_L\vc\left\{\Psi_{\bC}(\bs,\bxi_0)\right\}
+
b
\vc(\bSigma)
\left(
w_2(d_0)d_0^2-kw_3(d_0)
\right).
\end{split}
\]
Because $\Pi_L\vc(\bSigma)=\vc(\bSigma)$, we conclude that
\[
\bL\Lambda_{\btheta}'(\bxi_0)^{-1}
\Psi_{\btheta}(\bs,\bxi_0)
=
\Pi_L
\vc
\left\{
\Psi_{\bN}(\bs,\bxi_0)
\right\},
\]
where
\begin{equation}
\label{def:PsiN}
\Psi_{\bN}(\bs,\bxi)
=
v_1(d)(\by-\bX\bbeta)(\by-\bX\bbeta)^T
-
v_2(d)\bSigma,
\end{equation}
with $d^2=(\by-\bX\bbeta)^T\bV^{-1}(\by-\bX\bbeta)$
and
\begin{equation}
\label{def:v and w}
\begin{split}
v_1(s)
&=
aw_2(s)\\
v_2(s)
&=
-bw_2(s)s^2+(a+bk)w_3(s).
\end{split}
\end{equation}
Hence, if we define
\begin{equation}
\label{def:Nn}
\bN_n
=
\frac1n
\sum_{i=1}^{n}
\Psi_{\bN}(\bs_i,\bxi_0),
\end{equation}
with $\Psi_{\bN}$ defined in~\eqref{def:PsiN},
then it follows that
\[
\sqrt{n}(\vc(\bV_n)-\vc(\bSigma))
=
\Pi_L
\vc\left\{
\sqrt{n}(\bN_n-\E[\bN_n])
\right\}
+
o_P(1).
\]
This proves the first statement in Theorem~\ref{th:expansion}.

To prove the second statement,
note that from $\Lambda(\bxi_0)=\mathbf{0}$, 
together with~\eqref{eq:decomp psi-theta} and~\eqref{eq:kronecker},
it follows that
\[
\begin{split}
0
&=
\btheta_0^T
\E\left[\Psi_{\btheta}(\bs,\bxi_0)\right]
=
\E\left[
\vc(\bSigma^{-1})^T
\vc\left\{
\Psi_{\bC}(\bs,\bxi_0)
\right\}
\right]\\
&=
\E\left[
\text{tr}\left(\bSigma^{-1}\Psi_{\bC}(\bs,\bxi_0)\right)
\right]
=
\E\left[w_2(d_0)d_0^2-kw_3(d_0)\right],
\end{split}
\]
where $\Psi_{\bC}$ is defined in~\eqref{def:PsiC}.
Then, from the properties of elliptically contoured densities,
together with~\eqref{def:v and w}, one finds $\E[\Psi_{\bN}(\bs,\bxi_0)]=\mathbf{0}$.
This means that $\sqrt{n}(\bN_n-\E[\bN_n])$ is asymptotically normal with mean zero
and variance
\[
\E\left[
\vc(\Psi_{\bN}(\bs,\bxi_0))
\vc(\Psi_{\bN}(\bs,\bxi_0))^T
\right]
=
\E
\left[
\E\left[
\vc(\Psi_{\bN}(\bs,\bxi_0))
\vc(\Psi_{\bN}(\bs,\bxi_0))^T
\Big|\bX
\right]
\right].
\]
The inner expectation on the right hand side is the conditional expectation of $\by\mid\bX$,
which has the same distribution as $\bSigma^{1/2}\bz+\bmu$,
where $\bz$ has a spherical density $f_{\mathbf{0},\bI_k}(\bz)=g(\|\bz\|^2)$.
This implies that
\[
\begin{split}
&\E\left[
\vc\left(\bSigma^{-1/2}\Psi_{\bN}(\bs,\bxi_0)\bSigma^{-1/2}\right)
\vc\left(\bSigma^{-1/2}\Psi_{\bN}(\bs,\bxi_0)\bSigma^{-1/2}\right)^T
\right]\\
&=
\E_{\mathbf{0},\bI_k}\left[v_1(\|\bz\|)^2\|\bz\|^4\right]
\E_{\mathbf{0},\bI_k}\left[\vc\left(\bu\bu^T\right)\vc\left(\bu\bu^T\right)^T\right]\\
&\quad-
\E_{\mathbf{0},\bI_k}\left[v_1(\|\bz\|)v_2(\|\bz\|)\|\bz\|^2\right]
\E_{\mathbf{0},\bI_k}\left[\vc\left(\bu\bu^T\right)\vc\left(\bI_k\right)^T\right]\\
&\quad-
\E_{\mathbf{0},\bI_k}\left[v_1(\|\bz\|)v_2(\|\bz\|)\|\bz\|^2\right]
\E_{\mathbf{0},\bI_k}\left[\vc\left(\bI_k\right)\vc\left(\bu\bu^T\right)^T\right]\\
&\quad+
\E_{\mathbf{0},\bI_k}\left[v_2(\|\bz\|)^2\right]
\E_{\mathbf{0},\bI_k}\left[\vc\left(\bI_k\right)\vc\left(\bI_k\right)^T\right],
 \end{split}
\]
where $\bu=\bz/\|\bz\|$.
From Lemma~5.1 in~\cite{lopuhaa1989}, we have
\[
\mathbb{E}_{\mathbf{0},\bI_k}\vc(\bu\bu^T)\vc(\bu\bu^T)^T
=
\sigma_1(\bI_{k^2}+\mathbf{K}_{k,k})+\sigma_2\vc(\bI_k)\vc(\bI_k)^T,
\]
where $\sigma_1=\sigma_2=(k(k+2))^{-1}$.
Hence, the first term on the right hand side is equal to
\[
\frac{\E_{\mathbf{0},\bI_k}\left[v_1(\|\bz\|)^2\|\bz\|^4\right]}{k(k+2)}
\left(
\bI_{k^2}+\mathbf{K}_{k,k}+\vc(\bI_k)\vc(\bI_k)^T
\right).
\]
This leads to one term $\bI_{k^2}+\mathbf{K}_{k,k}$ with coefficient
\[
\frac{\E_{\mathbf{0},\bI_k}\left[v_1(\|\bz\|)^2\|\bz\|^4\right]}{k(k+2)}
\]
and using that, according to Lemma~11.4 in~\cite{supplement}, $\E_{\mathbf{0},\bI_k}\left[\bu\bu^T\right]=(1/k)\bI_k$,
we find a second term $\vc(\bI_k)\vc(\bI_k)^T$ with coefficient
\[
\frac{\E_{\mathbf{0},\bI_k}\left[v_1(\|\bz\|)^2\|\bz\|^4\right]}{k(k+2)}
-
\frac{2\E_{\mathbf{0},\bI_k}\left[v_1(\|\bz\|)v_2(\|\bz\|)\|\bz\|^2\right]}{k}
+
\E_{\mathbf{0},\bI_k}\left[v_2(\|\bz\|)^2\right].
\]
This means that
\[
\begin{split}
&\E\left[
\vc\left(\bSigma^{-1/2}\Psi_{\bN}(\bs,\bxi_0)\bSigma^{-1/2}\right)
\vc\left(\bSigma^{-1/2}\Psi_{\bN}(\bs,\bxi_0)\bSigma^{-1/2}\right)^T
\right]\\
&=
\sigma_1
\left(\bI_{k^2}+\mathbf{K}_{k,k}\right)
+
\sigma_2
\vc(\bI_k)\vc(\bI_k)^T
\end{split}
\]
where
\begin{equation}
\label{def:delta2}
\begin{split}
\sigma_1
&=
\frac{\E_{\mathbf{0},\bI_k}\left[v_1(\|\bz\|)^2\|\bz\|^4\right]}{k(k+2)}\\
\sigma_2
&=
\sigma_1
-
\frac{2\E_{\mathbf{0},\bI_k}\left[v_1(\|\bz\|)v_2(\|\bz\|)\|\bz\|^2\right]}{k}
+
\E_{\mathbf{0},\bI_k}\left[v_2(\|\bz\|)^2\right],
\end{split}
\end{equation}
or equivalently
\[
\E\left[
\vc(\Psi_{\bN}(\bs,\bxi_0))
\vc(\Psi_{\bN}(\bs,\bxi_0))^T
\right]
=
\sigma_1
\left(\bI_{k^2}+\mathbf{K}_{k,k}\right)(\bSigma^{-1}\otimes\bSigma^{-1})
+
\sigma_2
\vc(\bSigma)\vc(\bSigma)^T.
\]
We rewrite $\sigma_2$:
\[
\begin{split}
\sigma_2
&=
-\frac{2\sigma_1}{k}
+
\frac{(k+2)\sigma_1}{k}
-
\frac{2\E_{\mathbf{0},\bI_k}\left[v_1(\|\bz\|)v_2(\|\bz\|)\|\bz\|^2\right]}{k}
+
\E_{\mathbf{0},\bI_k}\left[v_2(\|\bz\|)^2\right]\\
&=
-\frac{2\sigma_1}{k}
+
\frac{\E_{\mathbf{0},\bI_k}\left[v_1(\|\bz\|)^2\|\bz\|^4\right]}{k^2}
-
\frac{2k\E_{\mathbf{0},\bI_k}\left[v_1(\|\bz\|)v_2(\|\bz\|)\|\bz\|^2\right]}{k^2}
+
\frac{k^2\E_{\mathbf{0},\bI_k}\left[v_2(\|\bz\|)^2\right]}{k^2}\\
&=
-\frac{2\sigma_1}{k}
+
\frac{1}{k^2}
\E_{\mathbf{0},\bI_k}\left[\Big(v_1(\|\bz\|)\|\bz\|^2-kv_2(\|\bz\|)\Big)^2\right].
\end{split}
\]
Furthermore
\[
v_1(s)s^2-kv_2(s)=(a+bk)\left\{w_2(s)s^2-kw_3(s)\right\},
\]
where
\[
a+kb
=
a+\frac{a\gamma_2}{\gamma_1-k\gamma_2}
=
\frac{1}{\gamma_1-k\gamma_2}
\]
where $\gamma_1$ and $\gamma_2$ are defined in~\eqref{def:gamma12}, which yields
\[
\gamma_1-k\gamma_2
=
\frac{1}{2k}
\E_{\mathbf{0},\bI_k}\Big[w_2'(\|\bz\|)\|\bz\|^3+2kw_3(\|\bz\|)-kw_3'(\|\bz\|)\|\bz\|\Big].
\]
It follows that
\[
\begin{split}
\sigma_1
&=
\frac{k(k+2)\E_{\mathbf{0},\bI_k}\left[w_2(\|\bz\|)^2\|\bz\|^4\right]}{
\Big(\E_{\mathbf{0},\bI_k}
\Big[
w_2'(\|\bz\|)\|\bz\|^3+k(k+2)w_3(\|\bz\|)
\Big]\Big)^2}\\
\sigma_2
&=
-\frac{2}{k}\sigma_1
+
\frac{4\E_{\mathbf{0},\bI_k}\left[\Big(w_2(\|\bz\|)\|\bz\|^2-kw_3(\|\bz\|)\Big)^2\right]}{
\left(\E_{\mathbf{0},\bI_k}\Big[w_2'(\|\bz\|)\|\bz\|^3+2kw_3(\|\bz\|)-kw_3'(\|\bz\|)\|\bz\|\Big]\right)^2.
}
\end{split}
\]
This proves the theorem.
\end{proof}

\paragraph{Proof of Theorem~\ref{th:limit H}.}
\begin{proof}
Let $H:\R^{k\times k}\to\R^{m}$ and let
\begin{equation}
\label{def:derivative H}
H'(\bV)=\frac{\partial H(\bV)}{\partial \vc(\bV)^T}=\left(\frac{\partial H_i(\bV)}{\partial v_{st}}\right)_{i=1,\ldots,m;\,s,t=1,\ldots,k}
\end{equation}
be the $m\times k^2$ matrix of partial derivatives.
According to the delta method
$\sqrt{n}(H(\bV_n)-H(\bSigma))$ is asymptotically normal with mean zero
and variance
$H'(\bSigma)
\text{\rm var}\{\text{\rm vec}(\bM)\}
H'(\bSigma)^T$.
Because $H$ is continuously differentiable and satisfies~\eqref{def:homogeneous mapping}, it follows that
\begin{equation}
\label{eq:derivative homogeneous mapping}
\sum_{j=1}^{l}
v_j
\frac{\partial H(\bv)}{\partial v_j}
=\mathbf{0}.
\end{equation}
This means that $H'(\bSigma)\vc(\bSigma)=\textbf{0}$.
Then, after inserting~\eqref{eq:limiting variance structured} for $\text{\rm var}\{\text{\rm vec}(\bM)\}$,
this finishes the proof of part~(i).

For part~(ii), let $H:\R^\ell\to\R^m$, and let
\begin{equation}
\label{def:derivative H theta}
H'(\btheta)=\frac{\partial H(\btheta)}{\partial \btheta^T}
=
\left(\frac{\partial H_i(\btheta)}{\partial \theta_j}\right)_{i=1,\ldots,m;\,j=1,\ldots,\ell}
\end{equation}
be the $m\times \ell$ matrix of partial derivatives.
According to the delta method
$\sqrt{n}(H(\btheta_n)-H(\btheta_0))$ is asymptotically normal with mean zero
and variance
\[
H'(\btheta_0)
\left\{
2\sigma_1\Big(\bL^T\left(\bSigma^{-1}\otimes\bSigma^{-1}\right)\bL\Big)^{-1}
+
\sigma_2
\btheta_0\btheta_0^T
\right\}
H'(\btheta_0)^T.
\]
Because $H$ satisfies~\eqref{def:homogeneous mapping} and~\eqref{eq:derivative homogeneous mapping},
it follows immediately that $H'(\btheta_0)\btheta_0=\textbf{0}$.
This finishes the proof of part~(ii).
\end{proof}

\paragraph{Proof of Lemma~\ref{lem:char IF}}
\begin{proof}
We apply Lemma~1 in~\cite{croux-haesbroeck2000}.
Although the lemma is established for the $N_k(\bmu,\bSigma)$ distribution,
the proof holds for any distribution with an elliptically contoured density.
According to~\cite{croux-haesbroeck2000}, there exist two functions $\alpha_C,\beta_C:[0,\infty)\to\R$,
such that
\begin{equation}
\label{eq:croux&haesbroeck}
\text{IF}(\by;\bC,P_{\bmu,\bSigma})
=
\alpha_C(d(\by))(\by-\bmu)(\by-\bmu)^T
-
\beta_C(d(\by))\bSigma.
\end{equation}
We have that
\[
\begin{split}
\text{\rm IF}(\by;\vc(\bM),P_{\bmu,\bSigma})
&=
\lim_{h\downarrow0}
\frac{\vc(\bM((1-h)P_{\bmu,\bSigma}+h\delta_{\by}))-\vc(\bM)(P_{\bmu,\bSigma})}{h}\\
&=
\Pi_L
\lim_{h\downarrow0}
\frac{\vc(\bC((1-h)P_{\bmu,\bSigma}+h\delta_{\by}))-\vc(\bC)(P_{\bmu,\bSigma})}{h}\\
&=
\Pi_L\vc\left(\text{IF}(\by;\bC,P_{\bmu,\bSigma})\right).
\end{split}
\]
Since $\bV$ is linear, it holds that $\vc(\bSigma)=\bL\btheta_0$ and because $\Pi_L$
is the projection on the column space of $\bL$, it follows that $\Pi_L\vc(\bSigma)=\vc(\bSigma)$.
When we insert the expression~\eqref{def:Projection matrix} for $\Pi_L$, together with~\eqref{eq:croux&haesbroeck}
and the fact that $(\bB^T\otimes\bA)\vc(\bv\bv^T)=\vc(\bA\bv\bv^T\bB)$ according to~\eqref{eq:kronecker}, this finishes the proof of part~(i).
Since $\bL$ has full column rank, $(\bL^T\bL)^{-1}\bL^T\vc(\bM(P))=\btheta(P)$, which yields
\[
\text{\rm IF}(\by;\btheta,P_{\bmu,\bSigma})
=
(\bL^T\bL)^{-1}\bL^T
\text{\rm IF}(\by;\vc(\bM),P_{\bmu,\bSigma}).
\]
Part~(i), together with~\eqref{eq:from vec to theta} finishes the proof of part~(ii).
\end{proof}

\paragraph{Proof of Lemma~\ref{lem:IF for H}}
\begin{proof}
Let $H:\R^{k\times k}\to\R^{m}$ with derivative $H'$ defined in~\eqref{def:derivative H}.
From the definition of influence function, it follows that
\begin{equation}
\label{eq:IF H(M)}
\begin{split}
\text{IF}(\by;H(\bM),P)
&=
\frac{\partial H(\bM(P_{h,\by}))}{\partial h}\bigg|_{h=0}\\
&=
\frac{\partial H(\bC)}{\partial \vc(\bC)^T}\bigg|_{\bC=\bM(P)}
\frac{\partial \vc(\bM(P_{h,\by}))}{\partial h}\bigg|_{h=0}\\
&=
H'(\bM(P))\cdot
\text{IF}(\by;\vc(\bM),P).
\end{split}
\end{equation}
Since $\bM(P)=\bSigma$, after inserting the expression in Lemma~\ref{lem:char IF},
together with $\vc(\bv\bv^T)=\bv\otimes\bv$,
for $\text{IF}(\by;\vc(\bM),P)$,
this proves part~(i).
Next, let $H:\R^{\ell}\to\R^m$ with derivative $H'$ defined by~\eqref{def:derivative H theta}.
It follows that
\begin{equation}
\label{eq:IF H(theta)}
\text{IF}(\by;H(\btheta),P)
=
H'(\btheta(P))
\cdot
\text{IF}(\by;\btheta,P).
\end{equation}
After inserting the expression in Lemma~\ref{lem:char IF}(ii) for $\text{IF}(\by;\btheta,P)$,
together with $\btheta(P)=\btheta_0$,
this proves part~(ii).
\end{proof}

\section{Derivation of $\sigma_1$ and $\sigma_2$}
\label{sec:comparison sigma}
We compare the expressions for $\sigma_1$ and $\sigma_2$ derived in Theorem~\ref{th:expansion}
with the ones obtained for specific cases in Tyler~\cite{tyler1982} and
Lopuha\"a \emph{et al}~\cite{lopuhaa-gares-ruizgazen2023}.

\paragraph{Example~\ref{example:LS}.}
Inserting $w_1=w_2=w_3=1$ in the expressions for $\sigma_1$ and $\sigma_2$ in Theorem~\ref{th:expansion}
gives
\[
\sigma_1
=
\frac{\E_{\mathbf{0},\bI_k}\left[\|\bz\|^4\right]}{k(k+2)},
\]
which equals 1 for the multivariate normal.
Furthermore,
\[
\begin{split}
\sigma_2
&=
-\frac{2}{k}
+
\frac{4\E_{\mathbf{0},\bI_k}\left[(\|\bz\|^2-k)^2\right]}{(2k)^2}\\
&=
-\frac{2}{k}
+
\frac{\E_{\mathbf{0},\bI_k}\left[\|\bz\|^4\right]-2k\E_{\mathbf{0},\bI_k}\left[\|\bz\|^2\right]+k^2}{k^2}\\
&=
-\frac{2}{k}
+
\frac{k(k+2)-2k^2+k^2}{k^2}
=0.
\end{split}
\]

\paragraph{Example~\ref{example:MLE}.}
First consider the special case of maximum likelihood, with
$w_1(s)=w_2(s)=-2g'(s^2)/g(s^2)$ and $w_3(s)=1$.
Note that
\begin{equation}
\label{eq:Lemma1 1997}
\E_{\mathbf{0},\bI_k}
\left[
z(\|\bz\|)
\right]
=
\frac{2\pi^{k/2}}{\Gamma(k/2)}
\int_{0}^{\infty}
z(r)
g(r^2)
r^{k-1}
\,\text{d}r,
\end{equation}
see e.g., Lemma~1 in Lopuha\"a~\cite{lopuhaa1997}.
When $\E_{\mathbf{0},\bI_k}[\|\bz^2\|]<\infty$,
then by means of integration by parts we get
\[
\begin{split}
\E_{\mathbf{0},\bI_k}
\Big[
w_2'(\|\bz\|)\|\bz\|^3
\Big]
&=
\frac{2\pi^{k/2}}{\Gamma(k/2)}
\int_{0}^{\infty}
\frac{4g'(r^2)^2}{g(r^2)^2}
g(r^2)
r^{k+3}
\text{d}r
-
\frac{2\pi^{k/2}}{\Gamma(k/2)}
\int_{0}^{\infty}
\frac{4g''(r^2)}{g(r^2)}
g(r^2)
r^{k+3}
\text{d}r\\
&=
4\E_{\mathbf{0},\bI_k}
\Big[
w_2(\|\bz\|)^2\|\bz\|^4
\Big]
-
k(k+2).
\end{split}
\]
It follows that
\begin{equation}
\label{eq:sigma1 tyler}
\sigma_1
=
\frac{k(k+2)\E_{\mathbf{0},\bI_k}\left[w_2(\|\bz\|)^2\|\bz\|^4\right]}{
\Big(\E_{\mathbf{0},\bI_k}[w_2'(\|\bz\|)\|\bz\|^3]
+
k(k+2)
\Big)^2}
=
\frac{k(k+2)}{\E_{\mathbf{0},\bI_k}\left[w_2(\|\bz\|)^2\|\bz\|^4\right]},
\end{equation}
which coincides with the expression found in Example~2 in Tyler~\cite{tyler1982},
who expresses expectations in terms of the random variable $T=\|\bz\|^2$.
To compute $\sigma_2$, first note that by means of integration by parts it follows
that $\E_{\mathbf{0},\bI_k}[w_2(\|\bz\|)\|\bz\|^2]=k$.
When we insert this in the expression for~$\sigma_2$,
this gives
\[
\begin{split}
\sigma_2
&=
-\frac{2}{k}\sigma_1
+
\frac{4\left(\E_{\mathbf{0},\bI_k}\left[w_2(\|\bz\|)^2\|\bz\|^4\right]-k^2\right)}{
\left(\E_{\mathbf{0},\bI_k}[w_2(\|\bz\|)^2\|\bz\|^4]-k(k+2)+2k\right)^2}
=
-\frac{2}{k}\sigma_1
+
\frac{4}{\E_{\mathbf{0},\bI_k}[w_2(\|\bz\|)^2\|\bz\|^4]-k^2}.
\end{split}
\]
After inserting $\E_{\mathbf{0},\bI_k}[w_2(\|\bz\|)^2\|\bz\|^4]=k(k+2)/\sigma_1$, 
as follows from~\eqref{eq:sigma1 tyler},
we find
\[
\sigma_2
=
-\frac{2}{k}\sigma_1
+
\frac{4}{k(k+2)/\sigma_1-k^2}
=
\frac{2\sigma_1(1-\sigma_1)}{k+2-k\sigma_1},
\]
which coincides with the expression found in Example~2 in Tyler~\cite{tyler1982}.

Next, consider the general case of M-estimators, with $w_3=1$.
First note that Tyler~\cite{tyler1982} uses a function $u_2$, which relates to our
function $w_2$ as $w_2(s)=u_2(s^2)$.
Then, since $\bxi_0=(\bbeta_0,\btheta_0)$ satisfies~\eqref{eq:score equations M for P},
we find that
\[
\begin{split}
0
&=
\btheta_0^T\E\left[\Psi_{\btheta}(\bs,\bxi_0)\right]\\
&=
\E\left[
\vc(\bSigma^{-1})^T
\vc
\left\{
w_2(d_0)
(\by-\bX\bbeta_0)(\by-\bX\bbeta_0)^T
-
\bSigma
\right\}
\right]\\
&=
\E\left[
\text{tr}
\left\{
w_2(d_0)
(\by-\bX\bbeta_0)(\by-\bX\bbeta_0)^T\bSigma^{-1}
-
\bI_k
\right\}
\right]\\
&=
\E_{\mathbf{0},\bI_k}
\left[
w_2(\|\bz\|)\|\bz\|^2
-
k
\right],
\end{split}
\]
where $d_0^2=(\by-\bX\bbeta_0)^T\bSigma^{-1}(\by-\bX\bbeta_0)$,
so that $k=\E_{\mathbf{0},\bI_k}[w_2(\|\bz\|)\|\bz\|^2]=\E[u_2(T)T]$.
It then follows that
\[
\begin{split}
\E_{\mathbf{0},\bI_k}\left[w_2(\|\bz\|)^2\|\bz\|^4\right]
&=
\E\left[u_2(T)^2T^2\right]
=
k(k+2)\psi_1\\
\E_{\mathbf{0},\bI_k}
\Big[
w_2'(\|\bz\|)\|\bz\|^3
\Big]
&=
2
\E\left[
u_2'(T)T^2
\right]
=
2k(\psi_2-1),
\end{split}
\]
where $\psi_1$ and $\psi_2$ are defined in Example~3 in Tyler~\cite{tyler1982}.
Then from the expressions provided in Theorem~\ref{th:expansion} we find
\[
\begin{split}
\sigma_1
&=
\frac{k^2(k+2)^2\psi_1}{
\left(
2k(\psi_2-1)+k(k+2)
\right)^2}
=
\frac{(k+2)^2\psi_1}{(2\psi_2+k)^2}\\
\sigma_2
&=
-\frac{2\sigma_1}{k}
+
\frac{4\left\{k(k+2)\psi_1-k^2\right\}}{(2k\psi_2)^2}.
\end{split}
\]
The expression for $\sigma_1$ coincides with the one in Example~3 in Tyler~\cite{tyler1982}.
After inserting this in~$\sigma_2$, one can verify that
also the expression for $\sigma_2$ coincides with one in Example~3 in Tyler~\cite{tyler1982}.

\paragraph{Example~\ref{example:S-estimators}.}
With $w_1(s)=\rho'(s)/s$, $w_2(s)=k\rho'(s)/s$ and $w_3(s)=\rho'(s)s-\rho(s)+b_0$,
one can easily verify that the expressions for $\sigma_1$ and $\sigma_2$ in Theorem~\ref{th:expansion}
coincide with the ones in Corollary~9.2 in Lopuha\"a \emph{et al}~\cite{lopuhaa-gares-ruizgazen2023}.

\section{Details for Examples~\ref{ex:shape} and~\ref{ex:direction}}
\label{sec:details shape direction}
\paragraph{Example~\ref{ex:shape}.}
From~\ref{eq:derivative shape H} we find
\begin{equation}
\label{eq:decomposition}
\begin{split}
&
H'(\bSigma)
\bL\Big(\bL^T\left(\bSigma^{-1}\otimes\bSigma^{-1}\right)\bL\Big)^{-1}\bL^T
H'(\bSigma)^T\\
&=
\frac{1}{k^2}|\bSigma|^{-2/k}
\vc(\bSigma)\vc(\bSigma^{-1})^T
\bL\Big(\bL^T\left(\bSigma^{-1}\otimes\bSigma^{-1}\right)\bL\Big)^{-1}\bL^T
\vc(\bSigma^{-1})\vc(\bSigma)^T\\
&\quad
-\frac{1}{k}|\bSigma|^{-2/k}
\vc(\bSigma)\vc(\bSigma^{-1})^T
\bL\Big(\bL^T\left(\bSigma^{-1}\otimes\bSigma^{-1}\right)\bL\Big)^{-1}\bL^T\\
&\quad
-\frac{1}{k}|\bSigma|^{-2/k}
\bL\Big(\bL^T\left(\bSigma^{-1}\otimes\bSigma^{-1}\right)\bL\Big)^{-1}\bL^T
\vc(\bSigma^{-1})\vc(\bSigma)^T\\
&\quad+
|\bSigma|^{-2/k}
\bL\Big(\bL^T\left(\bSigma^{-1}\otimes\bSigma^{-1}\right)\bL\Big)^{-1}\bL^T.
\end{split}
\end{equation}
Using that $\vc(\bSigma)=\bL\btheta_0$ and
\[
\vc(\bSigma^{-1})
=
\vc(\bSigma^{-1}\bSigma\bSigma^{-1})
=
\left(\bSigma^{-1}\otimes\bSigma^{-1}\right)
\vc(\bSigma)
=
\left(\bSigma^{-1}\otimes\bSigma^{-1}\right)
\bL\btheta_0,
\]
the first term on the right hand side of~\eqref{eq:decomposition} reduces to
$(1/k)
|\bSigma|^{-2/k}
\vc(\bSigma)
\vc(\bSigma)^T$.
Similarly, the second and third term on the right hand side
of~\eqref{eq:decomposition} are equal to
$-(1/k)
|\bSigma|^{-2/k}
\vc(\bSigma)
\vc(\bSigma)^T$.
Putting everything together, we find that the limiting covariance of $\sqrt{n}(H(\bV_n)-H(\bSigma))$
is given by~\eqref{eq:asymp variance shape}.

\paragraph{Example~\ref{ex:direction}.}
From Example~\ref{ex:shape} and the delta method, it follows that the limiting variance
of $\sqrt{n}(H(\btheta_n)-H(\btheta))$ is given by
\[
\begin{split}
&
(\bL^T\bL)^{-1}\bL^T
\left[
2\sigma_1
|\bSigma|^{-2/k}
\left\{
\bL\Big(\bL^T\left(\bSigma^{-1}\otimes\bSigma^{-1}\right)\bL\Big)^{-1}\bL^T
-
\frac1k
\vc(\bSigma)
\vc(\bSigma)^T
\right\}
\right]
\bL
(\bL^T\bL)^{-1}\\
&=
2\sigma_1
|\bSigma|^{-2/k}
\left\{
\Big(\bL^T\left(\bSigma^{-1}\otimes\bSigma^{-1}\right)\bL\Big)^{-1}
-
\frac1k
(\bL^T\bL)^{-1}\bL^T
\vc(\bSigma)
\vc(\bSigma)^T
\bL(\bL^T\bL)^{-1}
\right\}\\
&=
\frac{2\sigma_1}{|\bSigma|^{2/k}}
\left\{
\Big(\bL^T\left(\bSigma^{-1}\otimes\bSigma^{-1}\right)\bL\Big)^{-1}
-
\frac1k
\btheta_0\btheta_0^T
\right\},
\end{split}
\]
using that $\vc(\bSigma)=\bL\btheta_0$.


\begin{thebibliography}{10}

\bibitem{chervoneva2011}
I.~Chervoneva and M.~Vishnyakov.
\newblock Constrained {$S$}-estimators for linear mixed effects models with
  covariance components.
\newblock {\em Stat. Med.}, 30(14):1735--1750, 2011.

\bibitem{chervoneva2014}
I.~Chervoneva and M.~Vishnyakov.
\newblock Generalized s-estimators for linear mixed effects models.
\newblock {\em Statistica Sinica}, 24(3):1257--1276, 2014.

\bibitem{copt&heritier2007}
S.~Copt and S.~Heritier.
\newblock Robust alternatives to the f-test in mixed linear models based on
  mm-estimates.
\newblock {\em Biometrics}, 63(4):1045--1052, 2007.

\bibitem{copt2006high}
S.~Copt and M.~P. Victoria-Feser.
\newblock High-breakdown inference for mixed linear models.
\newblock {\em Journal of the American Statistical Association},
  101(473):292--300, 2006.

\bibitem{croux-haesbroeck2000}
C.~Croux and G.~Haesbroeck.
\newblock Principal component analysis based on robust estimators of the
  covariance or correlation matrix: influence functions and efficiencies.
\newblock {\em Biometrika}, 87(3):603--618, 2000.

\bibitem{fitzmaurice-laird-ware2011}
G.~M. Fitzmaurice, N.~M. Laird, and J.~H. Ware.
\newblock {\em Applied longitudinal analysis}.
\newblock Wiley Series in Probability and Statistics. John Wiley \& Sons, Inc.,
  Hoboken, NJ, second edition, 2011.

\bibitem{hampel1974}
F.~R. Hampel.
\newblock The influence curve and its role in robust estimation.
\newblock {\em J. Amer. Statist. Assoc.}, 69:383--393, 1974.

\bibitem{hartley&rao1967}
H.~O. Hartley and J.~N.~K. Rao.
\newblock Maximum-likelihood estimation for the mixed analysis of variance
  model.
\newblock {\em Biometrika}, 54:93--108, 1967.

\bibitem{huber1981}
P.~J. Huber.
\newblock {\em Robust statistics}.
\newblock Wiley Series in Probability and Mathematical Statistics. John Wiley
  \& Sons, Inc., New York, 1981.

\bibitem{jennrich&schluchter1986}
R.~I. Jennrich and M.~D. Schluchter.
\newblock Unbalanced repeated-measures models with structured covariance
  matrices.
\newblock {\em Biometrics}, 42(4):805--820, 1986.

\bibitem{kent&tyler1996}
J.~T. Kent and D.~E. Tyler.
\newblock Constrained {$M$}-estimation for multivariate location and scatter.
\newblock {\em Ann. Statist.}, 24(3):1346--1370, 1996.

\bibitem{kudraszow-maronna2011}
N.~L. Kudraszow and R.~A. Maronna.
\newblock Estimates of {MM} type for the multivariate linear model.
\newblock {\em J. Multivariate Anal.}, 102(9):1280--1292, 2011.

\bibitem{lopuhaa1989}
H.~P. Lopuha\"{a}.
\newblock On the relation between {$S$}-estimators and {$M$}-estimators of
  multivariate location and covariance.
\newblock {\em Ann. Statist.}, 17(4):1662--1683, 1989.

\bibitem{lopuhaa1997}
H.~P. Lopuha\"{a}.
\newblock Asymptotic expansion of {$S$}-estimators of location and covariance.
\newblock {\em Statist. Neerlandica}, 51(2):220--237, 1997.

\bibitem{lopuhaa2023}
H.~P. Lopuha\"a.
\newblock Highly efficient estimators with high breakdown point for linear
  models with structured covariance matrices.
\newblock {\em Econometrics and Statistics}, 2023.

\bibitem{lopuhaa-gares-ruizgazen2023}
H.~P. Lopuha\"a, V.~Gares, and A.~Ruiz-Gazen.
\newblock S-estimation in linear models with structured covariance matrices.
\newblock {\em Ann. Statist.}, 51(6):2415--2439, 2023.

\bibitem{supplement}
H.~P. Lopuha\"a, V.~Gares, and A.~Ruiz-Gazen.
\newblock Supplement to ``{S}-estimation in linear models with structured
  covariance matrices''.
\newblock 2023.

\bibitem{magnus&neudecker1988}
J.~R. Magnus and H.~Neudecker.
\newblock {\em Matrix differential calculus with applications in statistics and
  econometrics}.
\newblock Wiley Series in Probability and Mathematical Statistics: Applied
  Probability and Statistics. John Wiley \& Sons, Ltd., Chichester, 1988.

\bibitem{mallows1961}
C.~L. Mallows.
\newblock Latent vectors of random symmetric matrices.
\newblock {\em Biometrika}, 48:133--149, 1961.

\bibitem{mardia-marshall1984}
K.~V. Mardia and R.~J. Marshall.
\newblock Maximum likelihood estimation of models for residual covariance in
  spatial regression.
\newblock {\em Biometrika}, 71(1):135--146, 1984.

\bibitem{maronna1976}
R.~A. Maronna.
\newblock Robust {$M$}-estimators of multivariate location and scatter.
\newblock {\em Ann. Statist.}, 4(1):51--67, 1976.

\bibitem{miller1977}
J.~J. Miller.
\newblock Asymptotic properties of maximum likelihood estimates in the mixed
  model of the analysis of variance.
\newblock {\em Ann. Statist.}, 5(4):746--762, 1977.

\bibitem{rousseeuw-yohai1984}
P.~Rousseeuw and V.~Yohai.
\newblock Robust regression by means of {S}-estimators.
\newblock In {\em Robust and nonlinear time series analysis ({H}eidelberg,
  1983)}, volume~26 of {\em Lect. Notes Stat.}, pages 256--272. Springer, New
  York, 1984.

\bibitem{SalibianBarrera-VanAelst-Willems2006}
M.~Salibi\'{a}n-Barrera, S.~Van~Aelst, and G.~Willems.
\newblock Principal components analysis based on multivariate {MM} estimators
  with fast and robust bootstrap.
\newblock {\em J. Amer. Statist. Assoc.}, 101(475):1198--1211, 2006.

\bibitem{tatsuoka&tyler2000}
K.~S. Tatsuoka and D.~E. Tyler.
\newblock On the uniqueness of {$S$}-functionals and {$M$}-functionals under
  nonelliptical distributions.
\newblock {\em Ann. Statist.}, 28(4):1219--1243, 2000.

\bibitem{tyler1982}
D.~E. Tyler.
\newblock Radial estimates and the test for sphericity.
\newblock {\em Biometrika}, 69(2):429--436, 1982.

\bibitem{tyler1983}
D.~E. Tyler.
\newblock Robustness and efficiency properties of scatter matrices.
\newblock {\em Biometrika}, 70(2):411--420, 1983.

\bibitem{vanaelst&willems2005}
S.~Van~Aelst and G.~Willems.
\newblock Multivariate regression {$S$}-estimators for robust estimation and
  inference.
\newblock {\em Statist. Sinica}, 15(4):981--1001, 2005.

\end{thebibliography}

\end{document}